\documentclass[a4paper,11pt]{amsart}

\usepackage{epsfig}
\usepackage{color}
\usepackage{soul}
\usepackage{amssymb}

\newtheorem{theorem}{Theorem}
\newtheorem{corollary}[theorem]{Corollary}
\newtheorem{definition}{Definition}
\newtheorem{lemma}[theorem]{Lemma}
\newtheorem{proposition}[theorem]{Proposition} 
\newtheorem{remark}{Remark}

\newtheorem{maintheorem}{Theorem}

\newcommand{\CC}{{\mathbf C}}
\newcommand{\NN}{{\mathbf N}}
\newcommand{\ZZ}{{\mathbf Z}}
\newcommand{\RR}{{\mathbf R}}
\newcommand{\EU}{{\mathbf S}}

\newcommand{\dpt}{\displaystyle}
\usepackage[pagewise]{lineno}

\title[Three-dimensional horseshoes near an unfolding of a Hopf-Hopf singularity]{Three-dimensional horseshoes \\ near an unfolding of a Hopf-Hopf singularity}

\date{ \today}

\author[S. Ibanez and A. A. Rodrigues]{Santiago Ibanez$^{1}$ and Alexandre A. Rodrigues$^{1,2}$ \\
 \\ 
$^{1}$Departamento de Matem\'aticas, Universidad de Oviedo, \\ \medskip Av. Calvo Sotelo s/n, 33007 Oviedo, Spain \\
$^2$Lisbon School of Economics \& Management,\\
Centro de Matem\'atica Aplicada \`a Decis\~ao e Previs\~ao Econ\'omica \\Rua do Quelhas 6,  Lisboa 1200-781, Portugal 
 }

\thanks{The first author has been supported by the Spanish Research projects PID2020-113052GB-I00 and PID2023 147461NB-I00. The second author was supported by the Project CEMAPRE/REM -- UIDB /05069/2020 financed by FCT/MCTES through national funds. }
\email{ mesa@uniovi.es \quad alexandre.rodrigues@fc.up.pt, arodrigues@iseg.ulisboa.pt }
\begin{document}

\begin{abstract}

Motivated by a certain type of unfolding of a Hopf-Hopf singularity, we consider a one-parameter family $(f_\gamma)_{\gamma\geq0}$ of $C^3$--vector fields in $\RR^4$ whose flows exhibit a heteroclinic cycle associated to two  periodic solutions and a bifocus, all of them hyperbolic. It is formally proved that combining rotation with a generic condition concerning the transverse intersection between the three-dimensional invariant manifolds of the periodic solutions, all sets are highly distorted by the first return map and hyperbolic three-dimensional   horseshoes emerge, accumulating on the network. Infinitely many linked horseshoes prompt the coexistence of infinitely many saddle-type invariant sets for all values of $\gamma\gtrsim 0$  belonging to the heteroclinic class of the two hyperbolic periodic solutions.  We apply the results to a particular unfolding of the Hopf-Hopf singularity, the so called \emph{Gaspard-type unfolding}.   
  \end{abstract}

    \maketitle
  
  \bigbreak

\bigbreak
\textbf{2010 --- AMS Subject Classifications} 
{Primary:  34C37. Secondary: 34D20, 37C27}

\bigbreak

\section{Introduction}

The route through dynamical complexity started with Poincar\'e at the end of XIX century. In his seminal essay \cite{Poincare}, he discovered that homoclinic intersections between the invariant manifolds of a hyperbolic equilibrium were sources of  complex dynamics. In the thirties, Birkhoff \cite{Birkhoff} showed that, for   diffeomorphisms on the plane, there exists infinitely many periodic orbits near a transverse homoclinic intersection. In the sixties,  Smale \cite{Smale} conceived the famous horseshoe to explain the dynamics described by Birkhoff, via a conjugacy to a Bernoulli shift.  Mora and Viana \cite{MV93} proved the emergence of strange attractors when tangent homoclinic intersections are generically unfolded in families of diffeomorphisms on the plane. These non-hyperbolic and persistent (with respect to the Lebesgue measure)   attractors are ``like'' those described in \cite{BC91} for the H\'enon family.

On the other hand, Shilnikov \cite{Shilnikov65, Shilnikov67A, Shilnikov70} proved the existence of a countable set of periodic solutions (of saddle-type) in every small neighbourhood of a homoclinic cycle to a saddle-focus with positive divergence,  the counterpart of the Birkhoff results in the context of smooth vector fields on $\RR^3$. Charles Tresser \cite{Tresser} showed an infinite number of linked horseshoes   in every neighbourhood of this type of Shilnikov homoclinic cycle. When the homoclinic connection is broken (in a generic way), infinitely many horseshoes disappear, a phenomenon  accompanied by unfoldings of homoclinic tangencies to periodic orbits, leading to   persistent non-hyperbolic strange attractors as the ones described in \cite{MV93}.  In \cite{hom, PR_book, PR}, it was proven that infinitely many of these attractors can coexist for non generic families of vector fields. 
 Shilnikov  saddle-focus homoclinic orbits may be seen as  the simplest global configuration (in dimension 3) which can explain the existence of chaotic dynamics in families of vector fields; nonetheless their explicit occurrence   is not easy to prove analytically.
 
 Explicit results guaranteeing the existence of strange attractors are difficult to find in the literature. Although there are several examples based on numerical evidences, analytical proofs are rare. In the search for checkable and simple criteria to ensure the existence of observable chaos, we may find the \emph{unfoldings of singularities} of vector fields \cite{Barrientos_book}.  

 It has been proven that   saddle-focus homoclinic cycles appear in generic smooth unfoldings of certain singularities. Indeed, Ib\'a\~nez and Rodr\'iguez \cite{IbanezRodriguez} proved that these configurations have been found in generic  unfoldings of a three-dimensional nilpotent singularity of codimension three. This means that suspended H\'enon like strange attractors appear in generic unfoldings of such singularities. See   \cite{Barrientos} and references therein for additional details and \cite{dumibakok2005, dumibakok2001} for complementary results regarding the unfolding of the three-dimensional nilpotent singularity.  More recently, it has been proved the existence of Shilnikov saddle-focus homoclinic orbits in generic unfoldings of Hopf-zero type singularities of codimension two (see \cite{BIS, dumibakoksim}). In higher dimensions there  exist singularities of low codimension which can play the role of organizing centers of chaotic dynamics.
  For instance, the authors of  \cite{Drubi} provided numerical evidences of the existence of strange attractors in the unfolding of the codimension-three Hopf-Bogdanov-Takens bifurcation. This bifurcation has also been studied in the context of Hamiltonian systems.
  
One of the main motivations for obtaining simple criteria for the existence of strange attractors in given families of vector fields is their \emph{applicability.} Results in  Ib\'a\~nez and Rodr\'iguez \cite{IbanezRodriguez} have been  used to prove the existence of chaotic dynamics in a model of coupled oscillators  \cite{DIR07} and also in a general class of delay differential equations  \cite{CY08}.
 Family (1.2) in \cite{Barrientos}  has been widely studied in the literature because of its relevance in many physical settings as, for instance, the study of travelling waves in the Korteweg-de Vries model.\\

\subsection{Hopf-Hopf singularity}
In the present article, we are interested in  unfoldings of the Hopf-Hopf singularity whose $1$-jet at the origin of the coordinates $(x_1, y_1, x_2, y_2)\in \RR^4$, is linearly conjugate to
$$
\omega_1\left(x_1 \frac{\partial}{\partial y_1} -y_1 \frac{\partial}{\partial x_1}  \right)+ \omega_2\left(x_2 \frac{\partial}{\partial y_2} -y_2 \frac{\partial}{\partial x_2}  \right)
$$
with $\omega_1\neq \omega_2$ and $\omega_1, \omega_2 \in \RR^+$.

These singularities were classified by Takens \cite{Takens74} and, according to Kuznetsov \emph{et al} \cite{kuznetsov1998elements}, there is a topological type for which  the fifth order truncation of any generic unfolding may exhibit a heteroclinic configuration that could explain  chaotic dynamics (cf. \cite[\S 8.6.3]{kuznetsov1998elements}). For a suitable choice of parameters, the truncated system of order 5 has an hyperbolic bifocus $\textbf{O}$ and two hyperbolic non-trivial periodic solutions, say $\CC_1$ and $\CC_2$,  giving rise to an attracting heteroclinic network $\Gamma$. More precisely, in $\RR^4$, our starting point (organizing center) satisfies:

    \begin{equation*}
 \begin{array}{l}
 \dim W^u(\textbf{O}) = \dim W^s(\CC_2)=2,\\
  \dim W^u(\CC_2) = \dim W^s(\CC_1)=3,\\
    \dim W^u(\CC_1) = \dim W^s(  \textbf{O})=2\\
\end{array}
   \end{equation*}
and 
    \begin{equation*}
 \begin{array}{l}
 W^u(\bold{O}) \backslash \{\bold{O} \} \subset W^s(\CC_2),\\  W^u(\CC_2)\backslash \CC_2 \subset W^s(\CC_1) \\  W^s(\bold{O}) \backslash \{\bold{O} \} \subset W^u(\CC_1).
 \end{array}
   \end{equation*}

 The heteroclinic network $\Gamma \subset \RR^4$  can be seen as the organising center for the understanding of complex dynamics (in its unfoldings). Heteroclinic cycles and networks are flow-invariant sets that can occur robustly in dynamical systems and are often associated with intermittent behaviour \cite{ALR,   Rodrigues3, Rodrigues2022, RLA}.

 The dynamics near this kind of network is more intriguing that the ones studied in \cite{KR2010, Shilnikov65}  because of the existence of bifocus and the dimension 4.

The striking complexity of the dynamics near a   cycle containing a bifocus has been investigated by Shilnikov \cite{Shilnikov67A, Shilnikov70}, who claimed the existence of a countable set of periodic solutions of saddle type. It was shown that, for any $N \in \NN$ and for any local transverse section to the homoclinic cycle, there exists a compact invariant hyperbolic set on which the Poincar\'e map is topologically conjugate to the Bernoulli shift on $N$ symbols.  
In the works \cite{FS, LG}, the formation and bifurcations of periodic solutions were studied. The authors of \cite{IbRo} describe the hyperbolic suspended horseshoes that are contained in any small neighbourhood of a double homoclinic cycle to a bifocus and showed that \emph{heteroclinic switching} and \emph{suspended tridimensional horseshoes} are strongly connected. Examples of dynamical systems from applications where these homoclinic cycles play a basic role can be found in  \cite{Glendinning}. 
  In the reversible setting, the authors of \cite{Hart} proved the existence of a family of non-trivial (non-hyperbolic) closed trajectories and subsidiary connections near this type of cycle. See also the works by Lerman's team \cite{KL2009} (and references therein) who studied cycles to a bifocus in the hamiltonian context.  Results in \cite{Barrientos} show that homoclinic orbits to bifoci arise generically in unfoldings of four-dimensional nilpotent singularities of codimension 4.
  
 As already mentioned, one of the motivations for studying singularities is their applicability, particularly their role as organizing centers of complexity within a given model. In \cite{din25}, it is shown how Hopf-Hopf singularities can arise when two identical planar FitzHugh-Nagumo systems are coupled, leading to very interesting dynamical behaviors.  
  Compare with \cite{DIR11}, where the emergence of Hopf-pitchfork singularities was studied in the context of coupled oscillators -- see also \cite{druetal2024}. Hopf-Hopf singularities are also present in the unfolding of the four-dimensional nilpotent singularity of codimension four (see \cite{cdi2024}).

\subsection{Novelty and structure}
In this paper, motivated by certain unfoldings of a Hopf-Hopf singularity, we consider a one-parameter family $(f_\gamma)_{\gamma\geq0}$ of vector fields in $\RR^4$ whose flows exhibit a heteroclinic cycle associated to two  periodic solutions (with complex (non-real) Floquet multipliers) and a bifocus, all of them hyperbolic.

We assume that the flow of $\dot{x}= f_0(x)$ has an asymptotically stable (by ``inside'') heteroclinic cycle $\Gamma\equiv \Gamma_0$ associated to three hyperbolic saddles: two periodic solutions $\CC_1$, $\CC_2$ and a bifocus $\bold{O}$.
Their invariant manifolds coincide. For $\gamma\gtrsim 0$, the   connections involving  $\bold{O}$ persist while the others meet transversally into two disjoint curves. This gives rise to the network $\Gamma_\gamma$.

In this article, we state the existence of subsidiary heteroclinic connections from $\CC_2$ to $\CC_1$ (Proposition \ref{infinite number of connections}). Then, Theorem \ref{Th2} says that  the maximal invariant set associated to the first return map to a cross section to $\Gamma_\gamma$  has infinitely many  hyperbolic three-dimensional   horseshoes, accumulating on the network $\Gamma_\gamma$, all of them linked. This maximal invariant set could be seen as the  heteroclinic class of $\CC_2$ and $\CC_1$.
We apply the results to a particular unfolding of the Hopf-Hopf singularity, the so called \emph{Gaspard-type unfolding.}
In the present article, we prove partially the  claim of  \cite[pp.368]{kuznetsov1998elements}  that states the existence of Smale horseshoes near a generic Hopf-Hopf singularity.
  
Section \ref{s:HH} derives the computations performed by  \cite{kuznetsov1998elements} (the ``difficult case''). In Section \ref{description}, we describe rigorously the object and   hypotheses of our results, we state the main results and some consequences are discussed.
The $C^1$--normal forms are used in Section \ref{Local_dynamics} to construct three local transition maps around the saddles $\CC_1$, $\CC_2$ and the bifocus $\bold{O}$ and also global maps connecting the isolating blocks where the normal form is valid. Section \ref{s:global_geometry} deals with the geometrical structures which allow an understanding of the dynamics. We obtain results about the the geometry of  the return map and global consequences. This constitutes the main ingredients for the proof of the  Theorem \ref{Th2}. In Section  \ref{s:horseshoe}, we check precisely the Conley-Moser conditions to prove that, in the domain of definition of the first return map, there exists an infinite collection of hyperbolic horseshoes whose suspension accumulates on $\Gamma_\gamma$. Section   \ref{s:concluding} finishes the article.

 Throughout this paper, we have endeavoured to make a self contained exposition bringing together all topics related to the proofs. We have stated intermediate and understandable geometrical  lemmas and we have drawn illustrative figures to make the paper easily readable.

  \section{Motivation: Hopf-Hopf singularity \\ Derivation of the ``difficult case'' according to \cite{kuznetsov1998elements}}
   \label{s:HH}
In this section we revive the derivations of the unfoldings of a Hopf-Hopf singularity performed by Kuznetsov, with special focus on the parameters that yields a heteroclinic network, the so called ``difficult case'' \cite[pp. 361]{kuznetsov1998elements}.
Consider the system
  \begin{equation}
  \label{HH_general}
  \dot{x}= A(\alpha)x + F(x, \alpha)
  \end{equation}
  where $x=(x_1, x_2, x_3, x_4)\in \RR^4$, $\alpha=(\alpha_1, \alpha_2)\in \RR^2$ and $F(x, \alpha)=\mathcal{O}(\|x\|^2)$ is a $C^1$ map. 
   At $(\alpha_1, \alpha_2)=(0,0)$, we assume that $(x_1, x_2, x_3, x_4)=(0,0,0,0)$ is an equilibrium of \eqref{HH_general} such that, for $k\in \{1,2\}$, the linear map $A(\alpha)$ has two pairs of simple complex non-real eigenvalues of the form:
  $$
 \lambda_k(\alpha)= \mu_k (\alpha) \pm i \omega_k(\alpha). 
  $$
  such that $\mu_k(0)=0$ and $\omega_k(\alpha)>0$. The lowest codimension singularities in $\RR^4$ with a four-dimensional center manifold are the ones whose linear part is linearly conjugated to  
  $$
\left[ {\begin{array}{cccc}
   0& \omega_1 & 0& 0 \\
  -\omega_1&0 & 0& 0 \\
   0&0& \omega_2 & 0  \\
      0&0& 0& -\omega_2  \\
  \end{array} } \right].
    $$
  From now on, we are going to use bipolar coordinates given by:
     \begin{equation}
     \label{bipolar1}
  \left\{ 
\begin{array}{l}
x_1= r_1 \cos (\varphi_1) \\
x_2= r_1 \sin (\varphi_1) \\
x_3= r_2 \cos (\varphi_2) \\ 
x_4= r_2 \sin (\varphi_2). \\
\end{array}
\right.
  \end{equation}
  
  Under the non-degeneracy conditions \textbf{(HH.0)--(HH.5)} of Kuznetsov \cite{kuznetsov1998elements},   the system \eqref{HH_general} is locally smoothly equivalent (near $x=(0,0,0,0)$) to the system:
  \begin{equation}
  \label{HH_expansion}
  \left\{ 
\begin{array}{l}
\dot r_1=  r_1(\mu_1 + p_{11}(\mu)r_1^2  +p_{12}(\mu)r_2^2 + s_{1}(\mu)r_2^4)+ \Phi_1(r_1, r_2, \varphi_1, \varphi_2; \mu_1, \mu_2)  \\
\dot r_2=  r_2(\mu_2 + p_{21}(\mu)r_1^2  +p_{22}(\mu)r_2^2 + s_{2}(\mu)r_1^4)+ \Phi_2(r_1, r_2, \varphi_1, \varphi_2; \mu_1, \mu_2)  \\
\dot \varphi_1=  \omega_1(\mu) +\Psi_1(r_1, r_2, \varphi_1, \varphi_2; \mu_1, \mu_2)  \\
\dot \varphi_2=  \omega_2(\mu) +\Psi_2(r_1, r_2, \varphi_1, \varphi_2; \mu_1, \mu_2)  \\
\end{array}
\right.
  \end{equation}
    where $\Phi_k=\mathcal{O}((r_1^2+r_2^2)^3)$, $\Psi_k \equiv o(1)$, $k\in \{1,2\}$, are $2\pi$-periodic in $\varphi_1$ and $\varphi_2$ and the coefficientes $p_{jk}(0)$ and $s_k(0)$, $j,k=1,2$ could be explicitly computed (see \cite{kuznetsov1998elements}).
    The notation $\mathcal{O}$ and $o$ stands for the Landau notation.
   Truncating the radial coordinates of  \eqref{HH_expansion} to the order 5 and omitting the dependence on $\mu_1, \mu_2$ yields the following system:
    \begin{equation}
  \label{HH_truncation}
  \left\{ 
\begin{array}{l}
\dot r_1=  r_1(\mu_1 + p_{11} r_1^2  +p_{12}r_2^2 + s_{1}r_2^4)  \\
\dot r_2=  r_2(\mu_2 + p_{21} r_1^2  +p_{22}r_2^2 + s_{2}r_1^4) \\
\dot \varphi_1=  \omega_1  \\
\dot \varphi_2=  \omega_2.  \\
\end{array}
\right.
  \end{equation}
  
\begin{figure}
\begin{center}
\includegraphics[height=7.9cm]{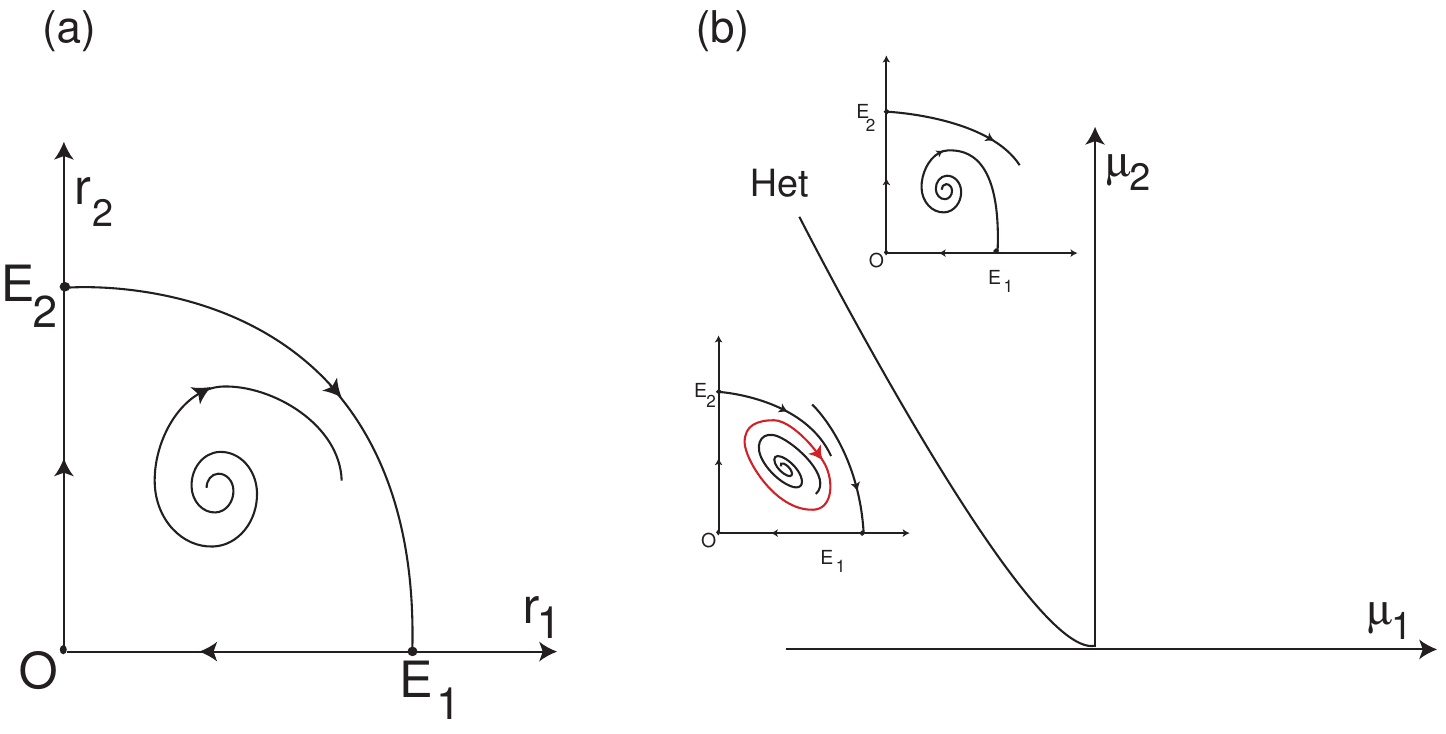}
\end{center}
\caption{\small  (a) Heteroclinic cycle associated to $\mathcal{O}\to (0,0)$ $\mathbf{E}_1\to (-\mu_1,0)$ and $ \mathbf{E}_2\to (0, \mu_2)$. The cycle is attracting by inside. (b) Bifurcation diagram of \eqref{HH_amplitude} in the parameters $(\mu_1,\mu_2)\in \RR^2$. The curve Het is the graph of the expression \eqref{curve1}.}
\label{Bif_diagram1}
\end{figure}

Ultimately one has to restore the tail which, generically, contains no $\varphi_1, \varphi_2$-dependent terms.
    Since the first two equations are independent of the last two, we are going to concentrate our attention to the truncated amplitude systems of order 5 associated to \eqref{HH_expansion}:
       \begin{equation}
  \label{HH_amplitude}
  \left\{ 
\begin{array}{l}
\dot r_1=  r_1(\mu_1 + p_{11} r_1^2  +p_{12}r_2^2 + s_{1}r_2^4)  \\
\dot r_2=  r_2(\mu_2 + p_{21} r_1^2  +p_{22}r_2^2 + s_{2}r_1^4)  . \\
\end{array}
\right.
  \end{equation}
  It is easy to check that the vector field associated to \eqref{HH_amplitude} is defined in $(\RR_0^+)^2$ and is $\ZZ_2(\gamma_1)\oplus \ZZ_2(\gamma_2)$-equivariant under the actions of 
  $$
  \gamma_1(r_1, r_2) = (-r_1, r_2)\quad \text{and}\quad \gamma_2(r_1, r_2) = (r_1, -r_2).
  $$
  In particular, the axes ${\rm Fix}(\gamma_1)$ and ${\rm Fix}(\gamma_2)$ are flow-invariant sets. If $p_{11}p_{22}<0$. (``difficult case'' of \cite{kuznetsov1998elements}), the flow of \eqref{HH_amplitude} has three saddle-equilibria depicted in Figure  \ref{Bif_diagram1}(a):
  
 $$
 \mathcal{O}\to (0,0), \quad  \mathbf{E}_1\to (-\mu_1,0),\quad \mathbf{E}_2\to (0, \mu_2),
 $$
  Finally, assuming that 
\begin{equation}
\label{condition1}
  \delta:=\frac{p_{21}}{p_{11}}<0, \qquad \theta:=\quad  \frac{p_{12}}{p_{22}}<0, \qquad  \delta.\theta:=\frac{p_{21}}{p_{11}} \frac{p_{12}}{p_{22}}>1
 \end{equation}
 and
 \begin{equation}
\label{condition2}
 (p_{21}(p_{21}-p_{11})s_1+p_{12}(p_{12}-p_{22})s_2)(0)<0,
 \end{equation}
 the authors of  \cite[cf. pp. 365] {kuznetsov1998elements} proved that, if $\dpt \Theta= \frac{s_1}{p_{22}^2}$, for parameter values along the curve
\begin{equation}
\label{curve1}
 {\rm Het}=\left\{(\mu_1, \mu_2)\in \RR^2: \mu_2= -\frac{\delta-1}{\theta-1}\mu_1-\frac{(\theta-1)\Delta+(\delta-1)\Theta}{(\theta-1)^3}\mu_1^2 +\mathcal{O}(\mu_1^3)\right\},
\end{equation}
  the
  unstable manifold of $ \mathbf{E}_2$ coincides with the stable manifold of $ \mathbf{E}_1$ -- see Figure \ref{Bif_diagram1}(b). This gives rise to an attracting heteroclinic cycle (attracting from ``inside'') associated to the three saddle-points $\mathbf{E}_1$, $\mathbf{E}_2$ and $\mathcal{O}$. 
  Coming back to the system \eqref{HH_truncation}, we may conclude that, along the curve  {\rm Het} (cf. \eqref{curve1}),  there exists a heteroclinic cycle associated to 
  the lift by two rotations  of $\mathbf{E}_1, \mathbf{E}_2$ and $\bold{O}$ that we call  $\CC_1$, $\CC_2$ and $\bold{O}$.
By construction, the resulting vector field, in rectangular coordinates, is $\textbf{SO}(2) \times \textbf{SO}(2)$--equivariant (details of the lift by rotation in \cite{RLA}).

 For $k\in \{1,2\}$, geometrically the solutions $\CC_k$ are $2\pi/\omega_k$ hyperbolic periodic solutions and $\bold{O}$ is a hyperbolic bifocus.

 Motivated by the derivation performed in this section, our starting point in the next section is a heteroclinic cycle in $\RR^4$ involving two hyperbolic periodic solutions and a bifocus. We state the existence of chaotic dynamics in a certain class of  unfolding of this configuration, pointing out that generic unfoldings of a Hopf-Hopf singularity might not belong to such class.

\section{Description of the problem and main results}
\label{description}
The object of our study is the dynamics around a general
heteroclinic network associated to a hyperbolic  bifocus and two non-trivial hyperbolic periodic solutions, for which we give a rigorous description here.
In the sequel, for $\gamma \geq 0$ (small),  we consider the  autonomous differential equation
\begin{equation}
\label{system1}
\dot{x}= f_\gamma(x)
\end{equation}
where $f_\gamma:\RR^4\to \RR^4$  is a $C^3$--vector field satisfying the  properties for $\gamma=0$ described in Subsection \ref{ss:oc}.
  We   use the notation $\overline{A}$ to denote the topological closure of a given set $A\subset\mathbb{R}^4$. \\

    \subsection{The organizing center ($\gamma=0$)}
    \,
    \label{ss:oc}
\medbreak
\begin{description}
 \item[\textbf{(P1)}]\label{P1}  The point $\bold{O}=(0,0,0,0)$ is an equilibrium point such that the spectrum of $df_0|_\bold{O}$ is $\{-C_0 \pm i\omega_1,E_0 \pm i\omega_2\}$ where  $C_0 > E_0 > 0$ and $\omega_1, \omega_2>0$.\\
 \medbreak
 \item[\textbf{(P2)}]\label{P2}
The set $\CC_1=\{(\cos(\omega_1 t), \sin (\omega_1 t), 0,0) \,:\, t\in [0, 2\pi/\omega_1)\}$ is  the parameterisation of a periodic solution of \eqref{system1} whose  non-trivial Floquet exponents  are given by $-C_1\pm i \omega_2$ and $E_1 $, where $C_1>E_1>0$.\\
 \medbreak
 \item[\textbf{(P3)}]\label{P3}
The set $\CC_2=\{(0,0,\cos(\omega_2 t), \sin (\omega_2 t)) \,:\, t\in [0, 2\pi/\omega_2)\}$ is the parameterisation of a periodic solution of \eqref{system1} whose  non-trivial Floquet exponents  are given by $-C_2 $ and $E_2\pm i \omega_1$, where $C_2>E_2>0$.\\
 \medbreak
 \item[\textbf{(P4)}]\label{P4} The following inclusions hold:\\
 $$W^u(\bold{O}) \backslash \{\bold{O} \} \subset W^s(\CC_2), \quad  W^u(\CC_2)\backslash \CC_2 \subset W^s(\CC_1)\backslash \CC_1, \quad \text{and} \quad W^s(\bold{O}) \backslash \{\bold{O} \} \subset W^u(\CC_1).$$
  \medbreak
  
    \end{description}

    \subsection{Effect of the perturbation governed by $\gamma$}
With respect to the effect of the parameter $\gamma$ on the dynamics,  we keep the notation $\bold O$, $\CC_1$ and $\CC_2$ to refer to the respective hyperbolic  continuations, and additionally assume that:\\
\begin{description}
 \item[\textbf{(P5)}]\label{P5}   The two-dimensional heteroclinic connections from $\bold{O}$ to $\CC_2$ and from $\CC_1$  to  $\bold{O}$ (cf. \textbf{(P4)})  persist.\\
 \end{description}

 For $\gamma\neq 0$, by Kupka-Smale Theorem, generically the invariant three-dimensional manifolds $W^u(\CC_2)$ and $W^s(\CC_1)$ are transverse (intersecting or not). Throughout this article, we assume that: \\
 \begin{description}
 \item[\textbf{(P6)}]\label{P6}  For $\gamma\neq 0$, $W^u(\CC_2) \pitchfork W^s(\CC_1)$ along two different tubular sets $\mathcal{T}_1$ and $\mathcal{T}_2$ joining $\CC_2$ and $\CC_1$, giving rise to two circles (in appropriate cross-sections to the network),  say $\ell_1$ and $\ell_2$, sometimes called by Yvone-Villarceau circles\footnote{If $A,B$ are submanifolds of $\RR^4$, the notation $A\pitchfork B$ means that $A\cap B\neq  \emptyset$ and the tangent spaces of $A$ and $B$, computed at all intersection points, generate $\RR^4$.}. \\     \end{description}
 
\begin{remark}
 
The global transition maps of Subsection \ref{global transition} will be defined taking into account the Hypotheses \textbf{(P5)} and \textbf{(P6)}.
Hypothesis \textbf{(P6)} corresponds to one of the expected unfoldings from the coincidence when $\gamma=0$) of the three-dimensional invariant manifolds of the hyperbolic  periodic solutions $\CC_1$ and $\CC_2$.  
 \end{remark}

From \textbf{(P6)}, observe that the sets $\mathcal{T}_1$ and $\mathcal{T}_2$  shrink to $\CC_1$ and $\CC_2$ as they approach the periodic solutions. 
For $\gamma \in [0, \varepsilon]$, where $\varepsilon \gtrsim 0$, the heteroclinic network associated to $\CC_1$, $\CC_2$ and $\bold{O}$ is denoted by $\Gamma_\gamma$ and we set up the following notation:\\
\begin{eqnarray*}
 \delta_0&=&\frac{C_0}{E_0}, \quad  \delta_1=\frac{C_1}{E_1},\quad  \delta_2=\frac{C_2}{E_2} \qquad \text{(saddle-value of each saddle)}\\\\
 \xi&=& 1+\frac{C_1}{E_0}+\frac{C_0C_1}{E_0E_2},\\ \\
  \delta&=& \delta_0\delta_1\delta_2  \qquad \qquad \qquad \qquad \quad \qquad \text{(saddle-value of the network)}. \\
 \end{eqnarray*}
 \bigbreak
We also assume a technical condition: \\
 \begin{description}
 \item[\textbf{(P7)}]\label{P7}  
  For $i\in \{0, 1, 2\}$, there exist $d_{1}^{(i)},d_{2}^{(i)}\in \RR^+$ such that for all $m,n\in \ZZ$, the following Diophantine inequality holds:
 $$
 |m C_i-nE_i|> d_{1}^{(i)} (|m|+|n|)^{-d_{2}^{(i)}}.
 $$
      \end{description}
      \bigbreak
      \bigbreak

 The combination of Hypotheses \textbf{(P1)--(P3)} and \textbf{(P7)} define an open set in the space of parameters $(C_0, E_0, C_1, E_1, C_2, E_2)$. 
Observe that  the dynamics close to the heteroclinic networks $\Gamma_0$ and $\Gamma_\gamma$ ($\gamma >0$), are qualitatively different. While in the first case, one observes regular dynamics  (zero topological entropy) near $\Gamma_0$, in the second the dynamics are characterised by the existence of chaos (conjugate to shift dynamics). 
 
 In the $C^1$--topology, small changes in $\gamma$ may affect the local dynamics near $\CC_1$, $\CC_2$ and $\bold{O}$. Since all saddles are hyperbolic and the dynamics is dominated by its linear part, without loss of generality, we omit this dependence.

\subsection{Main results}
 For $\gamma\gtrsim 0$, the two-dimensional heteroclinic connections given by  \textbf{(P6)} are called \emph{primary connections} (giving zero turns around $\Gamma_\gamma$). 
  The others, whose existence is guaranteed by   the next result, will be called \emph{subsidiary}   heteroclinic connections from $\CC_2$ to $\CC_1$.  
  
   \begin{proposition}
 \label{infinite number of connections}
 If $\gamma \gtrsim 0$, under  Hypotheses \textbf{(P1)--(P7)} on  the equation \eqref{system1}, there exist infinitely many two-dimensional heteroclinic connections from $\CC_2$ to $\CC_1$. 
 \end{proposition}
 
The proof of Proposition \ref{infinite number of connections} will be performed in Subsection \ref{ss:proof1}. For $\gamma \gtrsim 0$, let $\mathcal{V} $ be a sufficiently small neighbourhood of $\Gamma_\gamma$,  and let $\Sigma$ be a  three-dimensional local  cross section to $\Gamma_\gamma$.
 We prove that there exists a first return map $$\mathcal{R}_\gamma : S_1 \cup S_2 \subset \Sigma \rightarrow \Sigma$$ defined on a pair of disjoint sets $S_1$ and $S_2$ whose closure contains $\ell_1$ and $\ell_2$ (cf. \textbf{(P6)}), respectively.
Our main result, whose proof is performed in Subsection \ref{suspended horseshoes}, reads as follows:

\begin{maintheorem}
\label{Th2}
Under Hypotheses \textbf{(P1)--(P7)} on  the equation \eqref{system1}, for $\gamma\gtrsim 0$, there exists an infinite sequence of pairwise disjoint invariant sets $\Lambda_N \subset S_1 \cup S_2$, with $N\in\mathbb{N}$, accumulating on $\ell_1, \ell_2$, such that $\mathcal{R}_\gamma|_{\Lambda_N}$ is topologically conjugate to a full shift over an alphabet of two symbols $\{1,2\}$.  The maximal invariant set $$\Lambda_\gamma := \bigcup_{N\in \NN} \Lambda_N$$ is a Cantor set with zero Lebesgue measure.
\end{maintheorem}

The sequence of horseshoes  $\Lambda_\gamma$ may be seen as the heteroclinic class $\overline{W^u(\CC_2)\pitchfork W^s(\CC_1)} \cap \Sigma$ and corresponds exactly to the set of points whose orbits stay in  $\mathcal{V} $ for all $t\in \RR$ (one might shrink the neighbourhood $\mathcal{V}$, if necessary).
  The choice of the domains $S_1$ and $S_2$ will be specified in Section   \ref{s:horseshoe}. \\ \\

 \subsection{Dynamical consequence: Heteroclinic switching}\, 
For system \eqref{system1}, let $\mathcal{A}= \{\CC_1, \CC_2, \bold{O}\}$. 

\begin{definition}
\label{Def2}
If $k\in\NN$, a \emph{finite path of order $k$} on $\Gamma_\gamma$  is a sequence
     $
     \{\xi_{1},\ldots,\xi_{k}\}
     $
     of   heteroclinic connections in $\Gamma_\gamma$ such that $\omega(\xi_i)=\alpha(\xi_{i+1})$, for all $i\in \{1, ..., k-1\}$. We use the notation $\sigma^{k}=\{\xi_1, ...,\xi_k \}$ for this type of finite path.  An infinite path is an infinite sequence $\{\xi_i\}_{i=1}^\infty$, denoted by $\sigma^\infty$.
     \end{definition}

For $\gamma\gtrsim 0$, let   $\mathcal{V}$  be an arbitrarily  small  neighbourhood of the network $\Gamma_\gamma$ as above and let $W_j\subset \mathcal{V} $ be a neighbourhood of  $P_j\in \mathcal{A}$, $j \in\{1, \ldots, n\}$. For each heteroclinic connection $\xi_i$ in $\Gamma_\gamma$, consider a point $p_i\in \xi_i$ and a neighbourhood $V_i\subset \mathcal{V} $ of $p_i$. The collection of these neighbourhoods should be pairwise disjoint.  In the next definition,  adapted from \cite{ALR}, we denote by $\varphi_{\gamma}(t,x)$ the  flow associated to \eqref{system1}, for $\gamma \gtrsim 0$.

\begin{definition}
Given neighbourhoods as above, for $k\in \NN$, we say that the trajectory of a point $q\in \RR^4$ \emph{follows the finite path of order} $k$, say
$\sigma^k$ as before, if there exist two monotonically increasing sequences of times $(t_{j})_{j\in \{1,\ldots,k+1\}}$ and $(z_{j})_{j\in \{1,\ldots,k\}}$ such that for all $j \in \{1,\ldots,k\}$, we have $t_{j}<z_{j}<t_{j+1}$ and:\\

\begin{enumerate}
\item[(i)]
$\varphi_\gamma (t,q)\subset \mathcal{V}$ for all $t\in \, \, [ t_{1},t_{k+1}]$;
\item[(ii)]
$\varphi_\gamma (t_{j},q) \in W_j$ for all $j \in \{1,\ldots,k+1\}$ and $\varphi_\gamma (z_{j},q)\in V_{j}$ for all $j \in \{1,\ldots,k\}$;
\item[(iii)] for all $j=1,\ldots,k$ there exists a proper subinterval $I\subset \, \, (z_{j},z_{j+1})$ such that, given $t\in \, \, (z_{j},z_{j+1})$, $\varphi_\gamma(t,q)\in W_{j+1}$ if and only if $t\in I$.\\
\end{enumerate}
\end{definition}
The notion of a trajectory following an infinite path can be stated similarly; when we refer to points that follow a path, we mean that their trajectories do it.  Based on  \cite{ALR}, we set:

\begin{definition} There is:
\label{Def1}
\begin{enumerate}
\item[(i)]  \emph{finite switching} of order $k$ near $\Gamma_\gamma$ if  for each finite path of order $k$, say $\sigma^k$, and for each neighbourhood $\mathcal{V} $, there is a trajectory in $\mathcal{V}$ that follows $\sigma^k$ and
\item[(ii)]  \emph{infinite forward switching} (or simply \emph{switching}) near $\Gamma_\gamma$  by requiring that for each infinite path and for each neighbourhood $\mathcal{V} $, there is a trajectory in $\mathcal{V}$ that follows it.
\end{enumerate}
\end{definition}

Switching has been defined for positive time; this is why it is called by \emph{forward switching}. We may define analogously \emph{backward switching} by reversing the direction of the variable $t\in \RR$.
 
An infinite path on $\Gamma_\gamma$ can be considered as a pseudo-orbit of (\ref{system1}) with infinitely many discontinuities. Switching near $\Gamma_\gamma$ means that any pseudo-orbit in $\Gamma_\gamma$ can be realized. In Homburg and Knobloch \cite{HK}, using \emph{connectivity matrices}, the authors gave an equivalent definition of switching, emphasising the possibility of coding all trajectories that remain in a given neighbourhood of the network in both finite and infinite times.

\begin{corollary}
\label{corol_switching}
Under Hypotheses \textbf{(P1)--(P7)}  on the equation \eqref{system1}, the network $\Gamma_\gamma$ ($\gamma>0$) exhibits infinite forward and backward switching.
\end{corollary}

The proof of Corollary \ref{corol_switching} follows from the proof of Theorem \ref{Th2} since solutions that realize switching can be found within the suspended horseshoes that appear near $\Gamma_\gamma$.
  Theorem \ref{Th2} also implies the coding for all solutions remaining in $\mathcal{V} $ for all time.
One knows that the admissible order of saddles to follow $\Gamma_\gamma$ is given by 
$$
 \ldots \to \CC_1, \to  \bold{O} \to \CC_2 \to \CC_1\to  \bold{O} \to \CC_2\to \ldots
$$
The randomness for the coding emerges in the choice of the connections $\gamma_1$ or $\gamma_2$ as a consequence of \textbf{(P6)}; if a solution follows $\gamma_1$ let us code it by $1$, and $2$ otherwise. 
Therefore, a solution whose  trajectory lies in $\mathcal{V} $ for all time may be coded by the number associated to the primary heteroclinic connection.
Among the words coding the solutions, we may identify periodic solutions, as well as heteroclinic connections from $\CC_2$ to $\CC_1$.
In particular periodic  switching may be realized by periodic solutions.

 \subsection{Gaspard-type unfoldings of a Hopf-Hopf singularity}
 \label{ss:Gaspard-type}
 In this section, we show the existence of tridimensional horseshoes in particular analytic unfoldings of a Hopf-Hopf singularity.   From now on,  we will refer to a four-dimensional vector field $f^\star$ in $\RR^4$ as a Hopf-Hopf singularity if it satisfies the following conditions:
\begin{itemize}
\item $O \equiv (0,0,0,0)$ is an equilibrium of $f^\star$ and 
\item the spectrum of $df^\star|_O$ is $\{\pm i \omega_1, \pm i \omega_2\}$, with $\omega_1,\omega_2>0$.
\end{itemize}
 
 We introduce the following definition motivated by Equations (2.15) and (2.16) of Gaspard \cite{Gaspard}.
 \begin{definition}
 A \emph{Gaspard-type unfolding} of  a Hopf-Hopf singularity $f^\star$ associated to the truncation of order 5 \eqref{HH_truncation}  has the form (in bipolar coordinates \eqref{bipolar1}):
 \begin{equation}
  \left\{ 
\begin{array}{l}
\dot r_1=  r_1(\mu_1 + p_{11} r_1^2  +p_{12}r_2^2 + s_{1}r_2^4) +\mathcal{O}(\gamma) \mathcal{H}_1(r_1, r_2)\\
\dot r_2=  r_2(\mu_2 + p_{21} r_1^2  +p_{22}r_2^2 + s_{2}r_1^4) +\mathcal{O}(\gamma) \mathcal{H}_2(r_1, r_2)\\
\dot \varphi_1=  \omega_1  +\mathcal{O}(\gamma) \mathcal{H}_3(r_1, r_2)\\
\dot \varphi_2=  \omega_2  + \mathcal{O}(\gamma) \mathcal{H}_4(r_1, r_2),\\
\end{array}
\right.
  \end{equation}
where  $\gamma \geq 0$ is a small parameter and moreover:

\begin{enumerate}
\item  the coefficients $p_{11}, p_{12}, p_{21}, p_{22}, s_1, s_2$ satisfy conditions  \eqref{condition1} and  \eqref{condition2};
\item  the parameters $(\mu_1, \mu_2)\in \text{Het}$ (cf \eqref{curve1}) and 
\item  $\mathcal{H}_1, \mathcal{H}_2, \mathcal{H}_3, \mathcal{H}_4 $ are:
\begin{enumerate}
\item arbitrary analytic functions such that    $\mathcal{H}_i(r_1,r_2)= \mathcal{O}(\|(r_1, r_2)\|^6)$, for $i=1,2$, and $\mathcal{H}_i(r_1,r_2)= \mathcal{O}(\|(r_1, r_2)\|)$ for $i=3,4$;
\item  equivariant under the reflexions:
$$
\kappa_1(r_1, r_2)=(-r_1, r_2)\quad \text{and}\quad \kappa_2(r_1, r_2)=(r_1, -r_2).
$$
\end{enumerate}
\end{enumerate}
\end{definition}
 Note that the Hopf-Hopf singularity $f^*$ corresponds to $\mu_1=\mu_2=\gamma=0$.
\begin{corollary}
\label{corol1}
Let $f $ be a  Gaspard-type unfolding of a Hopf-Hopf singularity $f^\star$ associated to the truncation of order 5 (see \eqref{HH_truncation}). Then the differential equation $\dot{x}= f_0(x)$ satisfies Properties \textbf{(P1)--(P5)}. Moroever,  for $\gamma>0$, if \textbf{(P6)--(P7)} are satisfied, then the flow of $\dot{x}= f_\gamma(x)$ exhibits a tridimensional   hyperbolic horseshoe.
\end{corollary}

 Corollary \ref{corol1} is a consequence of Theorem \ref{Th2}. Truncation of order 5 is the lowest   truncation order that guarantees that the heteroclinic cycle associated to  $\mathbf{E}_1, \mathbf{E}_2, \bold{O}$ is asymptotically stable ``by inside''. Chaos persists even when the connections involving the bifocus are broken.   
Corollary \ref{corol1} cannot be applied directly to generic unfoldings of a Hopf-Hopf singularity because, in principle,  it might not be possible to express a generic Hopf-Hopf singularities as a  smooth perturbation of a vector field whose flow has an attracting heteroclinic cycle.

\section{Local and global  maps}
\label{Local_dynamics}
 Using local coordinates near the bifocus $\bold{O}$ and the periodic solutions $\CC_1$ and $\CC_2$, we will provide a construction of local and global transition maps. Finally, a return map around the heteroclinic network $\Gamma_\gamma$ will be defined, for $\gamma \gtrsim 0$.

 \subsection{Local map near the bifocus $\bold{O}$}
 One needs the  normal form that is used when studying the general saddle-focus case. This normal form has been constructed in Appendix A of \cite{Shilnikov et al}. 
Let $f_\gamma$ as in (\ref{system1})  and let $A$ and $B$ the  $(2\times 2)$-matrices as in \cite{Shilnikov et al} that depend on the parameter $\gamma$. It is clear that
$$
 A(0)= \begin{pmatrix} -C_0 & \omega_1 \\ \omega_1 &-C_0  \end{pmatrix} \qquad \text{and} \qquad B(0)= \begin{pmatrix} E_0 & \omega_2 \\ \omega_2 &E_0  \end{pmatrix}.
$$
The generalization of Bruno's theorem may be stated as:
\begin{theorem}[Shilnikov \emph{et al} \cite{Shilnikov et al}, adapted]
\label{normal form1}
There is a local $C^{3}$ transformation near the bifocus $O$ such that in  the new
coordinates $(x,y)=((x_1,x_2),(x_3, x_4)),$ the system casts as follows
\begin{equation}
\label{eq1}
\left\{ 
\begin{array}{l}
\dot x = A(\gamma) x + h(x,y, \gamma)x, \\ 
\dot y = B(\gamma) y + g(x,y, \gamma)y
\end{array}
\right.
\end{equation}
where: 
\begin{itemize}
\item  $A(\gamma)$ and $B(\gamma)$ are $(2\times 2)$-matrices; 
\medbreak
\item $h,g$ are $C^{2}$-smooth with respect to $(x,y)$, their first derivatives are $C^{1}$-smooth with respect to $(x,y, \gamma)$ and 
\medbreak
\item the following identities are valid for every $x=(x_1, x_2)$, $y=(y_1, y_2)$ and $\gamma\approx 0$:
$$
h(0,0, \gamma)=0,\; \quad g(0,0, \gamma)=0,\; \quad  h(x,0,\gamma)=0,\; \quad  g(0,y, \gamma)=0
$$
and 
$$
h(0,y, \gamma)=0,\; \quad g(x,0, \gamma)=0.
$$
\end{itemize}
\end{theorem}
In the rest of the paper, we omit the   terms $h$ and $g$  (which tend to zero as $\gamma$ vanishes, as well as their derivatives).  We are assuming that the neighbourhood $V_\bold{O}$ in which the flow can be written as the linear part near $\bold{O}$ is the cartesian product of two closed two-dimensional disks.
\medbreak

 Let us consider bipolar coordinates $(r_1, \phi_1,r_2, \phi_2) \in [0,\varepsilon] \times \RR \pmod{2\pi}\times [0,\varepsilon] \times \RR\pmod{2\pi}$ in $V_{O}$ such that
$$
  \left\{ 
\begin{array}{l}
x_1=r_1 \cos(\phi_1)  \\
x_2=r_1 \sin(\phi_1) \\
 x_3=r_2 \cos(\phi_2)  \\
x_4=r_2 \sin(\phi_2).  \\
\end{array}
\right.
 $$
In the previous bipolar coordinates the local invariant manifolds are given by
$$
W^s_{loc}(\bold{O}) \equiv \{ r_2=0\} \qquad \text{and} \qquad   W^u_{loc}(\bold{O}) \equiv \{ r_1=0\}
$$
and we can write system (\ref{eq1}) as
\begin{equation}
\label{equation1}
\dot{r}_1=-C_0 r_1, \qquad \dot{\phi}_1=\omega_1, \qquad \dot{r}_2=E_0 r_2 \qquad \text{and} \qquad  \dot{\phi}_2=\omega_2.
\end{equation}
Solving the above system explicitly we get
 
$$
  \left\{ 
\begin{array}{l}
r_1(t)=r_1(0) e^{-C_0 t}  \\
\phi_1(t)=\phi_1(0)+ \omega_1 t  \\
 r_2(t)=r_2(0) e^{E_0 t} \\
\phi_2(t)=\phi_2(0)+ \omega_2t.  \\
\end{array}
\right.
 $$
  
\bigbreak

\begin{figure}
\begin{center}
\includegraphics[height=5.6cm]{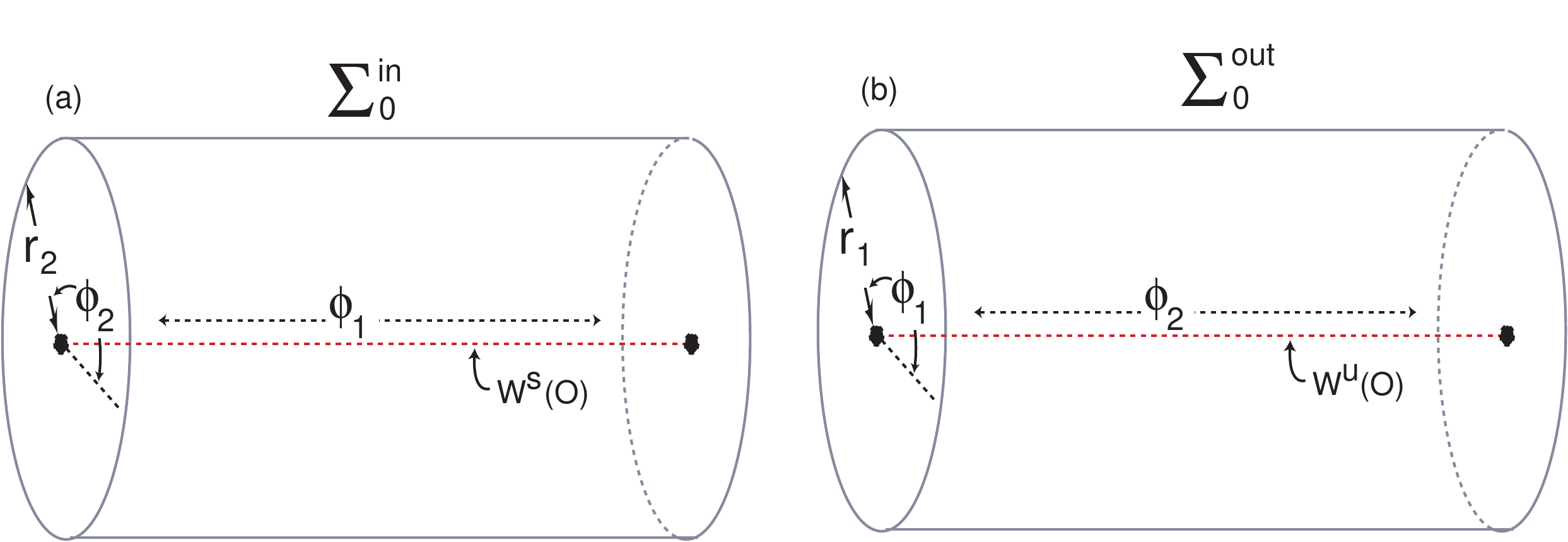}
\end{center}
\caption{\small  Coordinates and cross-sections near the bifocus $\bold{O}$. (a): $\Sigma_0^{in}$ and  (b): $\Sigma_0^{out}$. The superscripts ``in'' and ``out'' have been omitted in order to lighten the figures.}
\label{cross_sections1}
\end{figure}

Now, we consider two solid tori $\Sigma^{in}_0$ and $\Sigma^{out}_0$, transverse to the flow,  defined by:\\
\begin{enumerate}
\item[(a)] $\Sigma^{in}_0\equiv~\{r_1=\varepsilon \}$ with coordinates  $(\phi_1^{in}, r_2^{in}, \phi_2^{in})$ and
\medbreak
\item[(b)] $\Sigma^{out}_0\equiv \{r_2=\varepsilon \}$  with coordinates $(r_1^{out}, \phi_1^{out}, \phi_2^{out}).$\\
\end{enumerate}
 These two sets are illustrated in Figure \ref{cross_sections1}.
Solutions starting at interior points of $\Sigma^{in}_0$ go inside the hypertorus $V_\bold{O}$ in positive time; those starting at interior points of $\Sigma^{out}_0$ go outside $V_\bold{O}$ in positive time.
Intersections of local invariant manifolds at $\bold{O}$ and cross sections are circles parametrized as\\
\begin{equation}
\label{local s O}
W^s_{loc}(\bold{O}) \cap \Sigma^{in}_0=\{r_2^{in}=0,  0\leq \phi^{in}_1< 2\pi\}
\end{equation}
and
\begin{equation}
\label{local u O}
W^u_{loc}(\bold{O}) \cap \Sigma^{out}_0=\{r_1^{out}=0,  0\leq \phi^{out}_2 < 2\pi \}.
\end{equation}
 
The time of flight from $\Sigma^{in}_0$ to $\Sigma^{out}_0$ of a trajectory with initial condition $(\phi_1^{in}, r_2^{in}, \phi_2^{in}) \in \Sigma^{in}_O\backslash W^s_{loc}(\bold{O})$ is given by
$$
T=\frac{\ln(\varepsilon/ r_2^{in})}{E_0} = \frac{\ln \varepsilon -\ln r_2^{in} }{E_0}.
$$
Since $r_2^{in} \in (0,\varepsilon]$, $T$ is well defined and non-negative. Hence the local transition map
$$
\Pi_{0}: \Sigma^{in}_O \backslash W^s_{loc}(\bold{O}) \rightarrow \Sigma^{out}_O
$$
is given by
\begin{equation}
\label{local_flow_eq0}
\left(\begin{array}{l}\phi_1^{in}\\ \\r_2^{in}\\ \\\phi_2^{in}\end{array}\right)
\mapsto
\left(\begin{array}{l}r_1^{out}\\ \\\phi_1^{out} \\ \\ \phi_2^{out}\end{array}\right)
=
\left(\begin{array}{l}
\varepsilon^{1-\delta_0}\left(r_2^{in}\right)^{\delta_0}
\\
\\
\phi_1^{in}+\frac{\omega_1}{E_0} \ln (\varepsilon / r_2^{in}) \pmod{2\pi}
\\
\\
\phi_2^{in}+\frac{\omega_2}{E_0} \ln (\varepsilon / r_2^{in}) \pmod{2\pi}
\end{array}\right).
\end{equation}

\subsection{Local map near the periodic solution $\CC_1$}

By \textbf{(P2)}, the periodic solution $\CC_1$ is parameterised by $$ (\cos(\omega_1 t), \sin (\omega_1 t), 0,0), \quad t\in  [0, 2\pi/\omega_1).$$ 
We are assuming that the neighbourhood $V_1$ in which the flow can be linearised near $\CC_1$ is the product of an annulus and a disk.
Considering again bipolar coordinates   $$(r_1, \phi_1,r_2, \phi_2) \in [1-\varepsilon, 1+\varepsilon] \times \RR \pmod{2\pi}\times [0,\varepsilon] \times \RR\pmod{2\pi}$$ on $V_{1}$ such that
$$
x_1=r_1 \cos(\phi_1), \qquad x_2=r_1 \sin(\phi_1), \qquad x_3=r_2 \cos(\phi_2)\qquad \text{and} \qquad
x_4=r_2 \sin(\phi_2).
$$
 In bipolar coordinates, the local invariant manifolds are given by
$$
W^s_{loc}(\CC_1) \equiv \{ r_1=1\} \qquad \text{and} \qquad   W^u_{loc}(\CC_1) \equiv \{ r_2=0\}
$$
and, following the procedure of  \cite{ALR} and  the arguments of  \cite{WO2011} (Appendix A of \cite{LR2017} contains an argument for $C^2$--conjugacies),  there exists a $C^1$-change of coordinates such that system (\ref{system1}) may be written as (we are using \textbf{(P7)}):
\begin{equation}
\label{equation1}
\dot{r}_1= E_1 (r_1-1), \qquad \dot{\phi}_1=\omega_1, \qquad \dot{r}_2=-C_1 r_2 \qquad \text{and} \qquad  \dot{\phi}_2=\omega_2.
\end{equation}
Solving the above system explicitly we get
$$
\left\{
\begin{array}{l}
r_1(t)=1-(1-r_1(0)) e^{E_1 t}  \\
\phi_1(t)=\phi_1(0)+ \omega_1 t \\
 r_2(t)=r_2(0) e^{-C_1 t}\\
\phi_2(t)=\phi_2(0)+ \omega_2 t. \\
\end{array}
\right.
 $$

We consider two solid tori $\Sigma^{in}_1$ and $\Sigma^{out}_1$, transverse to the flow,  defined by:\
\begin{enumerate}
\item[(a)] $\Sigma^{in}_1\equiv~\{r_2=\varepsilon \}$ with coordinates  $( r_1^{in}, \phi_1^{in}, \phi_2^{in})$,
\medbreak
\item[(b)] $\Sigma^{out}_1\equiv \{r_1=1-\varepsilon \}$  with coordinates $(r_2^{out}, \phi_1^{out}, \phi_2^{out}).$
\end{enumerate}
 
\begin{figure}
\begin{center}
\includegraphics[height=6.4cm]{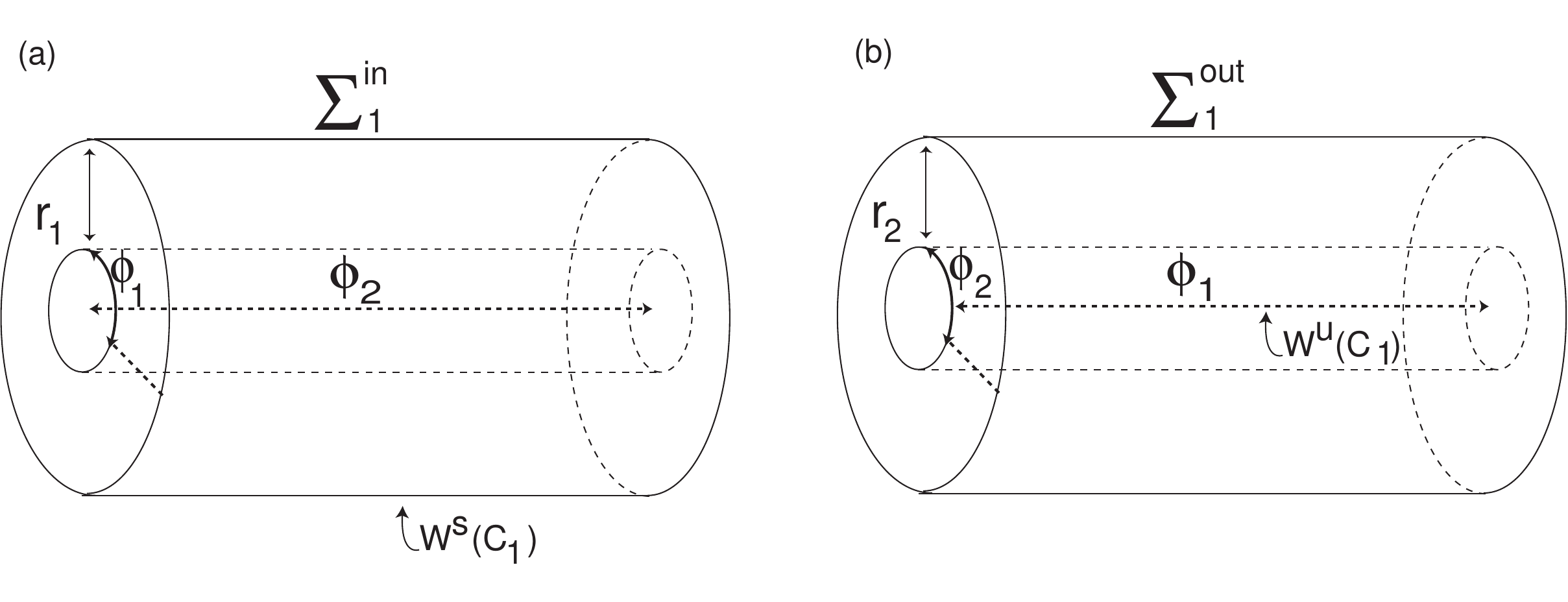}
\end{center}
\caption{\small  Coordinates and cross-sections near the periodic solution $\CC_1$. (a): $\Sigma_1^{in}$ and  (b): $\Sigma_1^{out}$. The superscripts ``in'' and ``out'' have been omitted in order to lighten the figures.}
\label{cross_sections2}
\end{figure}

\bigbreak
These two sets are sketched in Figure \ref{cross_sections2}.
By construction, trajectories starting at interior points of $\Sigma^{in}_1$ go inside  $V_1$ in positive time. Trajectories starting at interior points of $\Sigma^{out}_1$ go outside $V_1$ in positive time.
Intersections of local invariant manifolds at $\CC_1$ and cross sections are a two-dimensional torus parametrised as
\begin{equation}
\label{local s 1}
W^s_{loc}(\CC_1) \cap \Sigma^{in}_1=\{r_1^{in}=1, r_2^{in}=\varepsilon,  0\leq \phi^{in}_1, \phi^{in}_2< 2\pi\}
\end{equation}
and a circle given by:
\begin{equation}
\label{local u 1}
W^u_{loc}(\CC_1) \cap \Sigma^{out}_1=\{r_1^{out}=1-\varepsilon, r_2^{out}=0,  0\leq \phi^{out}_1, \phi^{out}_2< 2\pi \}.
\end{equation}

 From now on, we are interested in the case $ r_1^{in} \in [1-\varepsilon,1)$. 
The time of flight from $\Sigma^{in}_1\backslash W^s(\CC_1)$ to $\Sigma^{out}_1$ of a trajectory with initial condition $(r_1^{in}, \phi_1^{in},  \phi_2^{in}) \in \Sigma^{in}_1\backslash W^s_{loc}(\CC_1)$ is given by
$$
T=\frac{\ln \left(\frac{\varepsilon}{(1-r_1^{in})}\right)}{E_1} = \frac{\ln \varepsilon -\ln (1- r_1^{in})}{E_1}
$$
Since $r_1^{in} \in [1-\varepsilon,1)$, $T$ is well defined and non-negative. Hence the local transition map
$$
\Pi_{1}: \Sigma^{in}_1 \backslash W^s_{loc}(\CC_1) \rightarrow \Sigma^{out}_1
$$
is given by
\begin{equation}
\label{local_flow_eq1}
\left(\begin{array}{l}r_1^{in} \\\\\phi_1^{in}\\  \\\phi_2^{in}\end{array}\right)
\mapsto
\left(\begin{array}{l}r_2^{out}\\ \\\phi_1^{out} \\ \\ \phi_2^{out}\end{array}\right)
=
\left(\begin{array}{l}
   \varepsilon^{1-\delta_1}(1-r_1^{in})^{\delta_1}
\\
\\
\phi_1^{in}+\frac{\omega_1}{E_1} \ln \varepsilon - \frac{\omega_1}{E_1} \ln (1- r_1^{in}) \pmod{2\pi}
\\
\\
\phi_2^{in}+\frac{\omega_2}{E_1} \ln \varepsilon - \frac{\omega_2}{E_1} \ln (1- r_1^{in}) \pmod{2\pi}
\end{array}\right).
\end{equation}

\subsection{Local map near the periodic solution $\CC_2$}
The computations performed in this section are similar to the previous subsection. For the sake of completeness, we present all the details. 
 By \textbf{(P3)}, the periodic solution $\CC_2$ is parameterized by $$ (0,0, \cos(\omega_2 t), \sin (\omega_2 t)), \quad t\in  [0, 2\pi/\omega_2).$$ 
We are assuming that the neighbourhood $V_2$ in which the flow can be linearized near $\CC_2$ is the product of an annulus and a disk.
Considering again bipolar coordinates   $$(r_1, \phi_1,r_2, \phi_2) \in [0,\varepsilon] \times \RR \pmod{2\pi}\times [1-\varepsilon, 1]  \times \RR\pmod{2\pi}$$ on $V_{2}$ such that
$$
x_1=r_1 \cos(\phi_1), \qquad x_2=r_1 \sin(\phi_1), \qquad x_3=r_2 \cos(\phi_2)\qquad \text{and} \qquad
x_4=r_2 \sin(\phi_2).
$$
 In bipolar coordinates the local invariant manifolds are given by
$$
W^s_{loc}(\CC_2) \equiv \{ r_1=0\} \qquad \text{and} \qquad   W^u_{loc}(\CC_2) \equiv \{ r_2=1\}
$$
and we can write system (\ref{system1}) as
\begin{equation}
\label{equation1}
\dot{r}_1= E_2 r_1, \qquad \dot{\phi}_1=\omega_1, \qquad \dot{r}_2=-C_2 (r_2-1)\qquad \text{and} \qquad  \dot{\phi}_2=\omega_2.
\end{equation}
Solving the above system explicitly we get
$$
r_1(t)= r_1(0) e^{E_2 t}\qquad\phi_1(t)=\phi_1(0)+ \omega_1 t   \qquad  r_2(t)=1-(1-r_2(0)) e^{-C_2 t} \qquad   \phi_2(t)=\phi_2(0)+ \omega_2 t.
$$

We consider two solid tori $\Sigma^{in}_2$ and $\Sigma^{out}_2$ defined by
\begin{enumerate}
\item[(a)] $\Sigma^{in}_2\equiv~\{r_2=1-\varepsilon \}$ with coordinates  $( r_1^{in}, \phi_1^{in}, \phi_2^{in})$,
\medbreak
\item[(b)] $\Sigma^{out}_2\equiv \{r_1=\varepsilon \}$  with coordinates $(r_2^{out}, \phi_1^{out}, \phi_2^{out}).$
\end{enumerate}
 
\bigbreak
By construction, trajectories starting at interior points of $\Sigma^{in}_2$ go inside  $V_2$ in positive time. Trajectories starting at interior points of $\Sigma^{out}_2$ go outside $V_2$ in positive time.
Intersections of local invariant manifolds at $\CC_2$ and cross sections are  a circle parametrised as
\begin{equation}
\label{local s 1}
W^s_{loc}(\CC_2) \cap \Sigma^{in}_2=\{r_1^{in}=0, r_2^{in}=1-\varepsilon,  0\leq \phi^{in}_1, \phi^{in}_2< 2\pi\}
\end{equation}
and a two-dimensional torus parametrised by:
\begin{equation}
\label{local u 1}
W^u_{loc}(\CC_2) \cap \Sigma^{out}_2=\{r_1^{out}=\varepsilon, r_2^{out}=1,  0\leq \phi^{out}_1, \phi^{out}_2< 2\pi \}.
\end{equation}

 \bigbreak
The time of flight from $\Sigma^{in}_2$ to $\Sigma^{out}_2$ of a trajectory with initial condition $(r_1^{in}, \phi_1^{in},  \phi_2^{in}) \in \Sigma^{in}_2\backslash W^s_{loc}(\CC_2)$ is given by
$$
T=\frac{\ln\frac{\varepsilon}{r_1^{in}}}{E_2} = \frac{\ln \varepsilon -\ln  r_1^{in}}{E_2}.
$$
Since $r_1^{in} \in (0,\varepsilon]$, $T$ is well defined and non-negative. Hence the local transition map
$$
\Pi_{2}: \Sigma^{in}_2 \backslash W^s_{loc}(\CC_2) \rightarrow \Sigma^{out}_2
$$
is given by
\begin{equation}
\label{local_flow_eq2}
\left(\begin{array}{l}r_1^{in} \\\\\phi_1^{in}\\  \\\phi_2^{in}\end{array}\right)
\mapsto
\left(\begin{array}{l}r_2^{out}\\ \\\phi_1^{out} \\ \\ \phi_2^{out}\end{array}\right)
=
\left(\begin{array}{l}
 1-  \varepsilon^{1-\delta_2}(r_1^{in})^{\delta_2}
\\
\\
\phi_1^{in}+\frac{\omega_1}{E_2} \ln \varepsilon - \frac{\omega_1}{E_2} \ln ( r_1^{in}) \pmod{2\pi}
\\
\\
\phi_2^{in}+\frac{\omega_2}{E_2} \ln \varepsilon - \frac{\omega_2}{E_2} \ln ( r_1^{in}) \pmod{2\pi}
\end{array}\right).
\end{equation}

\subsection{Global transition maps and return maps}
\label{global transition}
Here, we define three global maps in disjoint neighbourhoods of the connections, corresponding to a flow-box around each connection.  We make use of the Hypotheses \textbf{(P5)} and \textbf{(P6)}.
For $\gamma=0$, we assume the following orientable global maps $$\Psi_{02}: \Sigma^{out}_0\to \Sigma^{in}_2, \quad \Psi_{21}: \Sigma^{out}_2\to \Sigma^{in}_1  \quad \text{and}\quad \Psi_{10}: \Sigma^{out}_1\to \Sigma^{in}_0$$ explicitly defined as: 
  $$
\begin{array}{lll}
(r_1^{out},  \phi_1^{out}, \varepsilon, \phi_2^{out})&\overset{\Psi_{02}}\mapsto& (A_1 r_1^{out}, \phi_1^{out}+A_2,1-\varepsilon,  \phi_2^{out}+A_3)\\ \\ 
(\varepsilon,  \phi_1^{out}, r_2^{out}, \phi_2^{out})&\overset{\Psi_{21}}\mapsto& (B_1 r_2^{out}, \phi_1^{out}+B_2, \varepsilon,   \phi_2^{out}+B_3) \\ \\
(1-\varepsilon, \phi_1^{out},  r_2^{out},  \phi_2^{out})&\overset{\Psi_{10}}\mapsto& (  \varepsilon,  \phi_1^{out}+C_1, C_2r_2^{out}, \phi_2^{out}+C_3) \\ \\
 \end{array}
 $$
 where $A_1, A_2, A_3, B_1, B_2, B_3, C_1, C_2, C_3,\in \RR^+$. For the sake of simplicity, we assume that $A_1=B_1=C_2=1$ and $A_2= A_3= B_2=B_3= C_1= C_3=0$.
  For $\gamma \gtrsim 0$, we assume that $\Psi_{02}$ and $\Psi_{10}$ remain the same and $\Psi_{21}$ can be written in cartesian rectangular coordinates as\\ 
 $$
 \left\{
\begin{array}{l}
X_1=D_1\varepsilon \cos(\phi_1^{out})  \\ \\
X_2=D_1\varepsilon \sin(\phi_1^{out})  \\\\
 X_3=D_2r_2^{out} \cos(\phi_2^{out})  \\ \\
 X_4=D_2r_2^{out} \sin(\phi_2^{out})  \\ 
\end{array}
 \right.
 $$
  where $D_1, D_2\in \RR^+$ and the transition map is given by:
$$
 (X_1, X_2, X_3, X_4) \mapsto ( \gamma +X_1, X_2, X_3, X_4). \\ \\
$$

In the proof of Theorem 1, we will not make explicit use of the expressions of the global transitions. This is why the previous assumptions are not real hypotheses of the problem.
 
\subsection{Notation and terminology for the case $\gamma \gtrsim 0$}
\label{ss:notation1}
\begin{figure}
\begin{center}
\includegraphics[trim={0 3cm 0 3cm},clip,width=13cm]{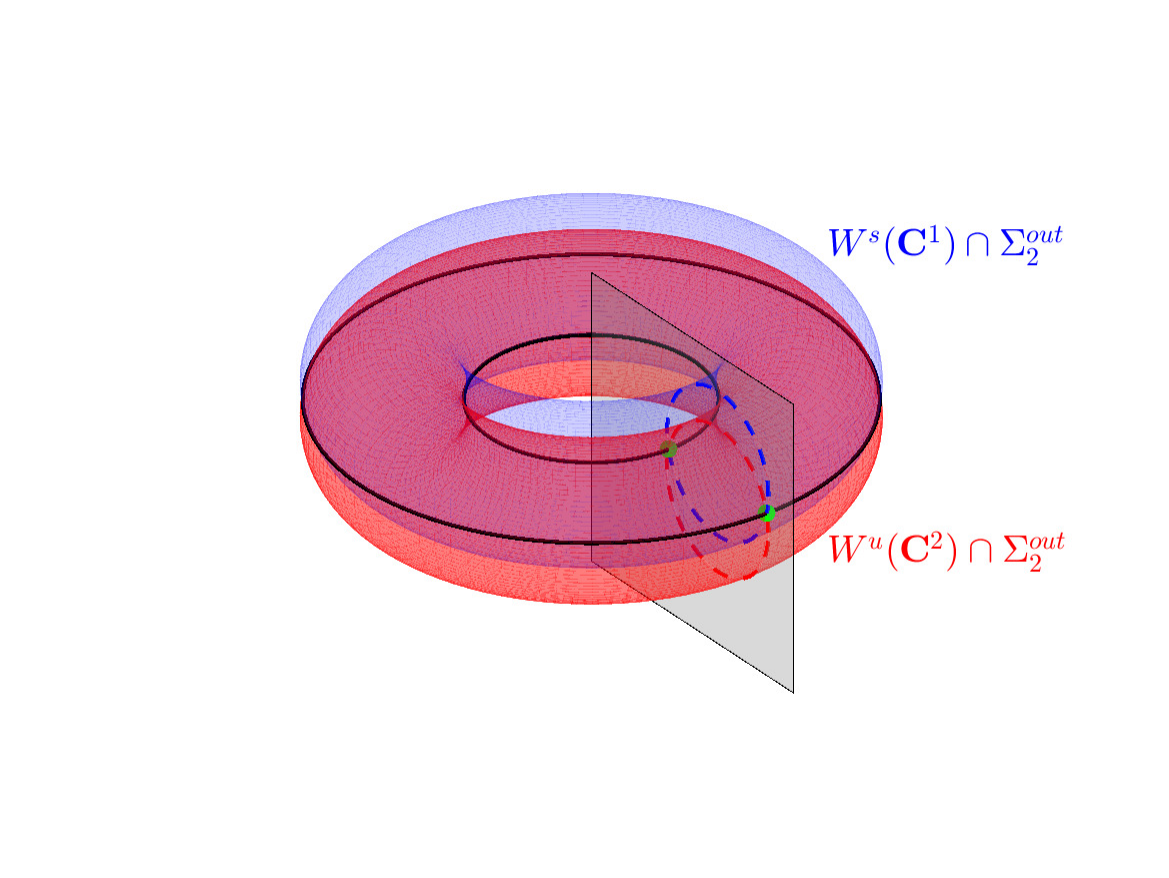}
\end{center}
\caption{\small Illustration of the intersections of $W^u(\mathcal{C}_2)$ and $W^s(\mathcal{C}_1)$ with $\Sigma_2^{out}$. The corresponding two-dimensional tori unfold from a coincidence, intersecting along two circles, $\ell_1^{out}$ and $\ell_2^{out}$. Extending these circles by the flow one generates the tubular sets $\mathcal{T}_1$ and $\mathcal{T}_2$. For greater clarity, a rectangle $P$, transverse to the tori, is included; the black circles in $P$ represent the intersections of the tori with $P$.}
\label{tori1}
\end{figure}

According to Hypothesis \textbf{(P6)}, we know that $W^u(\CC_2) \pitchfork W^s(\CC_1)$ along  two different tubular sets $\mathcal{T}_1$ and $\mathcal{T}_2$.
For $i \in \{1,2\}$, by taking $V_{1}$ and $V_2$ small enough, we can assume that $\mathcal{T}_i$ intersects each one of the cross sections $\Sigma_1^{in}$ and $\Sigma_2^{out}$ at exactly two closed curves (two circles), $\ell_i^{in}$ and $\ell_i^{out}$ defined by:
$$
 \ell_i^{in}=\mathcal{T}_i \cap \Sigma^{in}_1 \qquad \mbox{and} \qquad \ell_i^{out}=\mathcal{T}_i \cap \Sigma^{out}_2.
$$
 Figure \ref{tori1} provides an illustration of how the circles $\ell_1^{out}$ and $\ell_2^{out}$ emerge as the organizing center unfolds.

\begin{figure}
\begin{center}
\includegraphics[height=4.9cm]{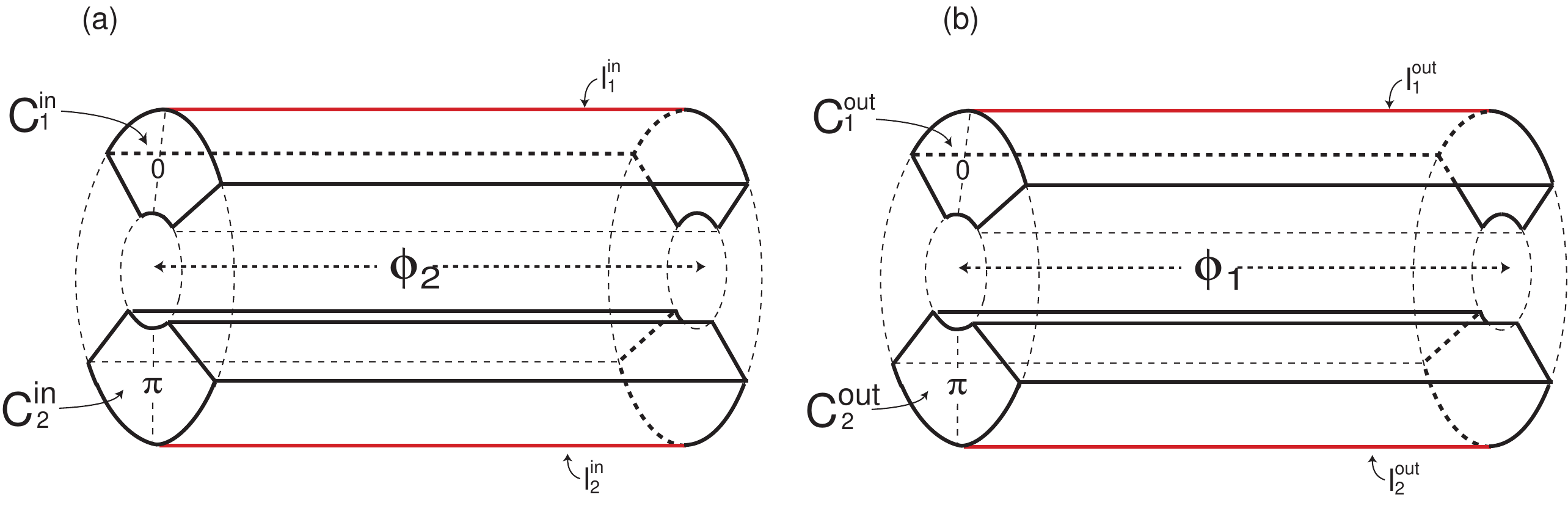}
\end{center}
\caption{\small  (a): Representation of $C_1^{in}$ and $C_2^{in}$ as subsets of $\Sigma_1^{in}$.  (b): Representation of $C_1^{out}$ and $C_2^{out}$ as subsets of $\Sigma_2^{out}$.}
\label{C1in}
\end{figure}

Let $\theta_i^{in}$ (resp. $\theta_i^{out}$) be the $\phi_1^{in}$ (resp. $\phi_2^{out}$) coordinate of $\ell_i^{in}$ (resp. $\ell_i^{out}$). For $i=1,2$, we can define the neighbourhoods  $C_i^{in}$ and $C_i^{out}$ as:
\begin{eqnarray*}
C_i^{in}=\left\{ (r_1^{in}, \phi_1^{in}, \phi_2^{in}):
\phi_1^{in} \in [\theta_i^{in}-\varepsilon^{in},\theta_i^{in}+\varepsilon^{in}],
r_1^{in}\in [1-\varepsilon^{in},1],
\phi_2^{in} \in [0,2\pi] \right\} \subset \Sigma_1^{in},
\end{eqnarray*}
for some small constant $1 \gg \varepsilon^{in}>0$, and \\
\begin{equation}
\label{ciout}
C_i^{out}=\left\{ ( \phi_1^{out}, r_2^{out}, \phi_2^{out}):
r_2^{out}\in [1-\varepsilon^{out},1],
\phi_1^{out} \in [0,2\pi], \phi_2^{out} \in [\theta_i^{out}-\varepsilon^{out},\theta_i^{out}+\varepsilon^{out}] \right\} \subset \Sigma_2^{out}
\end{equation}
for a small enough constant $1\gg \varepsilon^{out}>0$. The sets $C_1^{in}$ and $C_2^{in}$ are illustrated in Figure  \ref{C1in}. \\

 For the sake of simplicity, we assume that $\theta_1^{out}=0=\theta_1^{in}$ and $\theta_2^{out}=\pi=\theta_2^{in}$. 
From now on, let us denote the map $$\Pi_2\circ \Psi_{02}\circ \Pi_0\circ \Psi_{10}\circ \Pi_1: \Sigma_1^{in} \to \Sigma_2^{out}$$ by $\mathcal{G}$, which does not depend on the parameter $\gamma$. If necessary, we might need to shrink the domain and the range in order that the map is well defined. For $i\in \{1,2\}$, define the sets 
\begin{equation}
\label{def_S1}
S_i= \mathcal{G}^{-1}(C_i^{out}\backslash W^u(\CC_2))\subset \Sigma_1^{in}.
\end{equation}
Define the transition map $\Psi_{21}^{(i)}: C_i^{out} \to C_i^{in}$ (given the set $C_i^{in}$,   the tridimensional set $C_i^{out}$ exists and it is well defined by the Flow-box Theorem -- cf. \cite{pm}).
If $\varepsilon^{in}, \varepsilon^{out}>0$ are small enough, then $$C_1^{in}\cap C_2^{in}=\varnothing, \quad C_1^{out}\cap C_2^{out}=\varnothing \quad \text{and}\quad S_1\cap S_2=\varnothing.$$
For $\gamma \gtrsim 0$, we define the map $\Psi_{12}: C_1^{out}\cup C_2^{out} \to \Sigma_1^{in}$ as follows:
$$
\Psi_{12} ( \phi_1^{out}, r_2^{out}, \phi_2^{out})= 
\left\{\begin{array}{ll}
\Psi_{12}^{(1)} (   \phi_1^{out}, r_2^{out}, \phi_2^{out}) & \text{if} \quad   ( \phi_1^{out}, r_2^{out}, \phi_2^{out}) \in C_1^{out} \\ \\
\Psi_{12}^{(2)} (   \phi_1^{out}, r_2^{out}, \phi_2^{out}) & \text{if} \quad   ( \phi_1^{out}, r_2^{out}, \phi_2^{out}) \in C_2^{out}  \\
 \end{array} 
 \right.
 .$$
 Note that, for $\gamma\gtrsim 0$, the domain of  $\Psi_{12}$ has two connected components.                                                          

\subsection{Useful results and notation}
\label{ss:notation2}
In what follows, we introduce the folllowing notation:
$$\xi = \frac{1}{E_1}\left(1+\frac{C_1}{E_0}+\frac{C_0C_1}{E_0E_2}\right) \quad \text{and} \quad k(\varepsilon)= \varepsilon^{1-\delta}.$$
Observe that
$$
\frac{\xi}{\delta}= \frac{1}{C_2}\left(1+\frac{E_2}{C_0}+\frac{E_0 E_2}{C_0 C_1}\right).
$$

 \begin{lemma}
 \label{lemma1}
For $\gamma \gtrsim 0$ fixed and $\varepsilon > 0$ small, the expressions of \\$$\mathcal{G}: \Sigma_1^{in} \backslash W^s(\CC_1)\to \Sigma_2^{out} \quad \text{and} \quad\mathcal{G}^{-1}: \Sigma_2^{out} \backslash W^u(\CC_2) \to \Sigma_1^{in}$$ \\do not depend on $\gamma$ and are explicitly given, respectively, by:

\begin{equation}
\label{return1}
\left(\begin{array}{l}r_1^{in} \\\\\phi_1^{in}\\  \\\phi_2^{in}\end{array}\right)
\mapsto
\left(\begin{array}{l}\phi_1^{out}\\ \\ r_2^{out} \\ \\ \phi_2^{out}\end{array}\right)
=
\left(\begin{array}{l}
    \phi_1^{in}+ \xi \omega_1 \ln (\varepsilon)- \xi \omega_1  \ln (1-r_1^{in})  \pmod{2\pi}
\\
\\
1-k (\varepsilon) (1- r_1^{in})^\delta
\\
\\
 \phi_2^{in}+ \xi \omega_2 \ln (\varepsilon)- \xi \omega_2 \ln (1-r_1^{in})  \pmod{2\pi}
\end{array}\right)
\end{equation}

and
  
  \begin{equation}
\label{return2}
\left(\begin{array}{l}\phi_1^{out} \\\\r_2^{out}\\  \\\phi_2^{out}\end{array}\right)
\mapsto
\left(\begin{array}{l}r_1^{in} \\\\\phi_1^{in}\\  \\\phi_2^{in}\end{array}\right)
=
\left(\begin{array}{l}
 1-   \left(\frac{1-r_2^{out}}{k(\varepsilon)} \right)^\frac{1}{\delta}
\\
\\
\phi_1^{out} -\frac{\xi \omega_1}{\delta}\ln \varepsilon +\frac{\xi \omega_1}{\delta}\ln (1-r_2^{out}) \pmod{2\pi}
\\
\\
\phi_2^{out} -\frac{\xi \omega_2}{\delta}\ln \varepsilon +\frac{\xi \omega_2}{\delta}\ln (1-r_2^{out}) \pmod{2\pi}
\end{array}\right).
\end{equation}
 
\end{lemma}

\bigbreak

\begin{proof}
The independence of $\mathcal{G}$ and $\mathcal{G}^{-1}$ on $\gamma$ is a consequence of  \textbf{(P5)}, \textbf{(P6)} and the way the transitions maps $\Psi_{10}$, $\Psi_{02}$ are defined. The proof follows from the composition map $\Pi_2\circ \Pi_0\circ \Pi_1$ where the explicit expressions of $\Pi_1$,  $\Pi_0$ and $\Pi_2$ are given by \eqref{local_flow_eq1}, \eqref{local_flow_eq0} and \eqref{local_flow_eq2}, respectively.

\end{proof}

 From now on, for the sake of simplicity we assume that $\varepsilon=1$ and, in particular, $\ln \varepsilon=0$ and $k(\varepsilon)=1$.
 The return map $\mathcal{R}_\gamma$ to $\Sigma_1^{in}\backslash W^s(\CC_1)$ is defined by $\Psi_{21}\circ \mathcal{G}$, restricted to the set $\tilde\Lambda_\gamma$. This set should be understood as the maximal subset of $\Sigma_1^{in}\backslash W^s(\CC_1)$ where the return is well defined.   We will see that, in fact, $\tilde\Lambda_\gamma=\Lambda_\gamma$, with $\Lambda_\gamma$ as given in Theorem \ref{Th2}.

\section{Global geometry}
\label{s:global_geometry}
 In this section, we introduce some geometrical objects that will be useful in the sequel as the main ingredients to prove Theorem \ref{Th2}.

\subsection{Spiralling behaviour --  geometric preliminaries}

\begin{definition}
\label{spiralling1_def}
Let $a,b\in \RR$  such that $a < b $ and let $H$ be an annulus parametrised by a cover $$(\theta,h) \in \RR\times [a,b]$$ where $\theta$ is a periodic covering. A spiral on $H$ accumulating on the circle $\mathcal{C}\subset H$ parameterized by $h = h_0 \in [a,b]$ is a smooth curve
$$\alpha :[c, +\infty[ \rightarrow  H,$$
with $c\in \RR$ and such that its coordinates $\alpha (s)=( \theta(s), h(s))$ are functions of $s$ and:
 \begin{enumerate}

\item the map $\theta$ is monotonic for some unbounded subinterval of $[c,+\infty[$ and
\item $\dpt \lim_{s\to +\infty}|\theta(s)|=+\infty$.
\item the map $h$ is bounded by two monotonically decreasing maps converging to $h_0$ as $s \rightarrow +\infty$. \\
\end{enumerate}
\end{definition}

\begin{figure}
\begin{center}
\includegraphics[height=4.1cm]{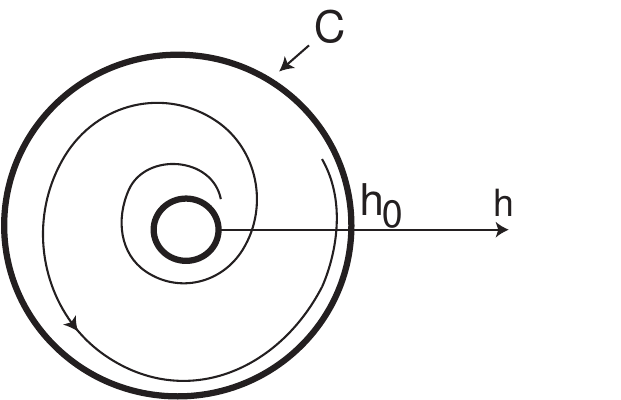}
\end{center}
\caption{\small Spiral on the annulus accumulating on the circle $\mathcal{C}$. The arrow indicates the direction of evolution of the variable $s$ from Definition \ref{spiralling1_def}.}
\label{spiral1}
\end{figure}

 Figure \ref{spiral1} illustrates this notion of a spiral accumulating on a circle. In the next definition we refer again to the circle $\mathcal{C}$ on the cylinder $H$ parametrised by $h=h_0\in (a,b)$ as introduced in Definition \ref{spiralling1_def}.\\
 
\begin{definition}
Let $H$ be an  annulus as in Definition \ref{spiralling1_def} and $c_1, c_2\in \RR$.
A \emph{double spiral} on $H$ around the curve $\mathcal{C}$ is the union of two spirals
$${ \alpha_i :[c_i, +\infty[ \rightarrow  H,}$$
  with $i=1,2$, accumulating on $\mathcal{C}$ and a curve connecting $\alpha_1(c_1)$ and $\alpha_2(c_2)$. The region bounded by a double spiral is called a \emph{spiralling region}. \\
\end{definition}

\begin{figure}
\begin{center}
\includegraphics[height=5.8cm]{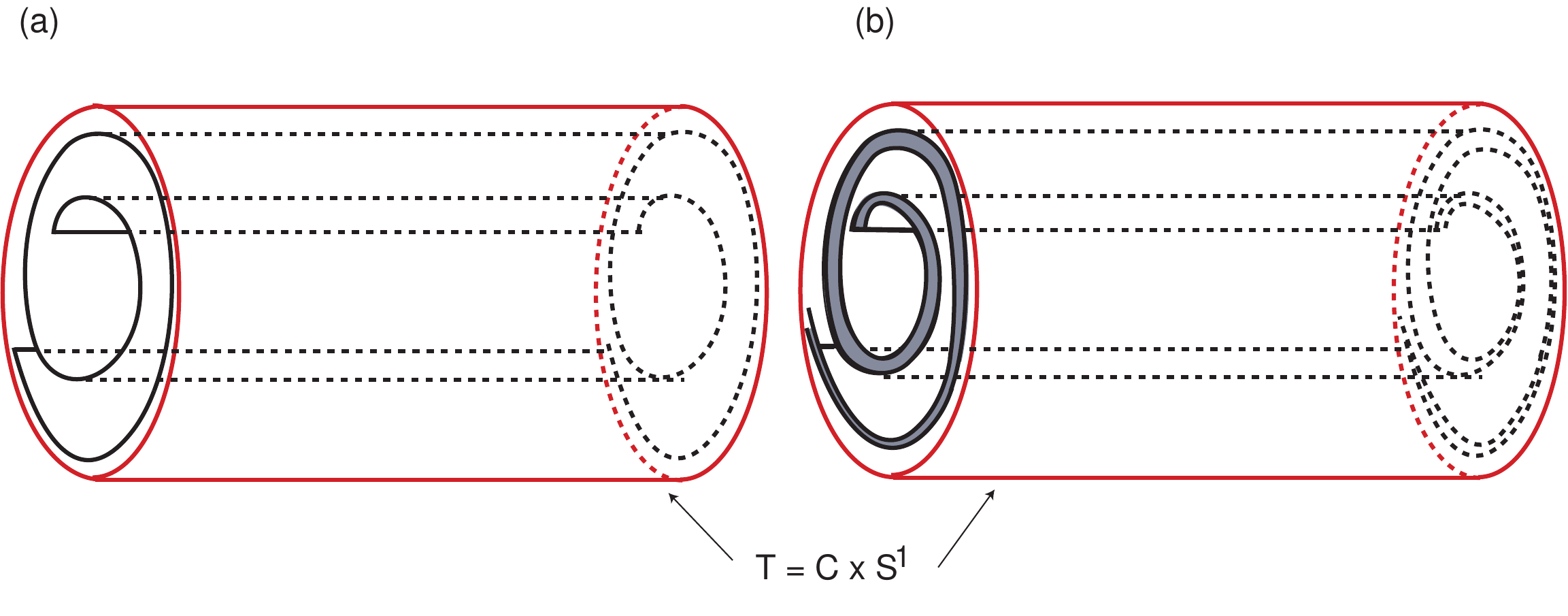}
\end{center}
\caption{\small (a) Spiralling sheet and (b) scroll, both accumulating on the torus $\mathcal{T}\equiv \mathcal{C}\times \EU^1$.}
\label{spiralling_sheet1}
\end{figure}

We introduce now the notion of \emph{spiralling sheet} depicted in Figure \ref{spiralling_sheet1}(a). \\

\begin{definition}
\label{ss_def2}
A two-dimensional manifold  $\mathcal{H}$ embedded in $\mathbb{R}^3$ is called a \emph{spiralling sheet} accumulating on a torus $\mathcal{T}\equiv \mathcal{C}\times \EU^1$  if there exist a spiral $S$ around $\mathcal{C}$, a neighbourhood $V\subset \mathbb{R}^3$ of $\mathcal{T}$, a neighbourhood $W_0 \subset \mathbb{R}^2$ of $\mathcal{C}$, a non-degenerated closed interval $I$ and a diffeomorphism $$\eta: V\rightarrow I \times W_0$$ such that:
$$\eta(\mathcal{H}\cap V)=I \times (S \cap W_0) \qquad \text{and} \qquad \mathcal{T}=\eta^{-1}(I\times \mathcal{C}).$$
\end{definition}

\bigbreak

According to Definition \ref{ss_def2}, up to a diffeomorphism, we may think on a \emph{spiralling sheet} accumulating on the torus  $\mathcal{T}\equiv \mathcal{C}\times \EU^1$ as the cartesian product of a spiral accumulating on $\mathcal{C}$ and $\EU^1$.  Observe that the diffeomorphic image of a spiralling sheet accumulating on a torus is again a spiralling sheet accumulating on another torus. In the present article, the torus $\mathcal{T}$ corresponds either to  $W^s(\CC_1)\cap \Sigma_1^{in}$ or $W^u(\CC_2)\cap \Sigma_2^{out}$, according to the scenario under analysis.\\
 
 \begin{remark}
Let $\mathcal{H}\subset \Sigma_2^{out}$ be a spiralling sheet accumulating in  $W^u_{loc}(\CC_2)\cap \Sigma_2^{out}$ and let    $P$ be a plane with  $\phi_1^{out}$ constant   cutting transversely  $W^u_{loc}(\CC_2)\cap \Sigma_2^{out}$. Then  the set $\mathcal{H}\cap P$ is a spiral accumulating on the circle $W^u_{loc}(\CC_2)\cap P$. The same holds for any spiralling sheet accumulating in  $W^s_{loc}(\CC_1)\cap \Sigma_1^{in}$.\\
 \end{remark}
 
\begin{definition}
\label{scroll_def}
Given two spiralling sheets $\mathcal{H}_1$ and $\mathcal{H}_2$ accumulating on a torus $\mathcal{T}= \mathcal{C}\times \EU^1\subset \RR^3$, any region limited by $\mathcal{H}_1$ and $\mathcal{H}_2$ inside a neighbourhood of $\mathcal{T}$ is said a \emph{scroll} accumulating on $\mathcal{T}$.
\end{definition}
Roughly speaking, as suggested by Figure \ref{spiralling_sheet1}(b), a scroll is a thick spiralling sheet, (i.e. ``with volume''). 

\subsection{Local geometry}
\label{firstreturn}
The following result shows that a two-dimensional annulus, or a set diffeomorphic to it, contained in $\Sigma_1^{in}$ and transverse to $W^{s}_{loc}(\CC_1) \cap \Sigma_1^{in}$ is sent by the map $\mathcal{G}$  into a spiralling sheet on $\Sigma_{2}^{out}$ accumulating on $W^u(\CC_2)\cap \Sigma_2^{out}$. An analogous result holds for  $\mathcal{G}^{-1}$.

\begin{proposition}
\label{My_Prop5}
For $\hat\delta>0$ arbitrarily small, let $\Xi: D\subset \mathbb{R}^2 \rightarrow \RR$ be a $C^1$ map defined on:
$$D=\{(u,v)\in\mathbb{R}^2 : 1-\hat\delta\leq u^2+v^2 \leq 1  \}$$
  which is constant for  $(u,v)\in D$ such that $u^2+v^2=1$, and let
$$
\mathcal{F}^{in}=\{(r_1^{in},\phi_1^{in},   \phi_2^{in}) \in \Sigma^{in}_1 : \phi_1^{in}=\Xi(r_1^{in} \cos \phi_2^{in},r_1^{in} \sin \phi_2^{in}),\,1-\hat\delta \leq r_1^{in}\leq 1,\,0\leq \phi_2^{in} <2\pi\},
$$
and 
$$
\mathcal{F}^{out}=\{(\phi_1^{out}, r_2^{out}, \phi_2^{out}) \in \Sigma^{out}_2 : \phi_2^{out}=\Xi(r_2^{out} \cos \phi_1^{out},r_2^{out} \sin \phi_1^{out}),\,1-\hat\delta \leq r_2^{out}\leq 1,\,0\leq \phi_1^{out} <2\pi\}. \\
$$
\bigbreak
Then the following assertions are valid:\\
\begin{enumerate}
\item  The set $\mathcal{G}(\mathcal{F}^{in}\backslash W^s_{loc}(\CC_1))$ is a spiralling sheet accumulating on $W^u_{loc}(\CC_2) \cap \Sigma^{out}_2$.\\
\item  The set $\mathcal{G}^{-1}(\mathcal{F}^{out}\backslash W^u_{loc}(\CC_2))$ is a spiralling sheet accumulating on $W^s_{loc}(\CC_1) \cap \Sigma^{out}_1$.\\
\item The sets $S_1$ and $S_2$ defined in \eqref{def_S1} are scrolls contained in $\Sigma_1^{in}$ accumulating on $W^{s}_{loc} (\CC_1)\cap \Sigma_1^{in}$. \\
 
\end{enumerate}
\end{proposition}

\begin{proof}
To simplify expressions, along this proof we use the following functions:
\begin{eqnarray}
\nonumber M(r_2^{out}) & = & 1-(1-r_2^{out})^\frac{1}{\delta}\\
\nonumber \mathcal{W}\equiv \mathcal{W}  (r_2^{out}, \phi_2^{out}) &=& \phi_2^{out} +\frac{\xi\omega_2}{\delta} \ln(1- r_2^{out})\\
\nonumber u\equiv u (r_2^{out}, \phi_2^{out})&=&  M(r_2^{out}) \cos\left(\mathcal{W}  (r_2^{out}, \phi_2^{out})\right) \\
\nonumber v\equiv v  (r_2^{out}, \phi_2^{out})&=&  M(r_2^{out}) \sin\left(\mathcal{W}  (r_2^{out}, \phi_2^{out})\right).
\nonumber
\end{eqnarray}
We perform the proof of the lemma item by item.

\bigbreak

\textit{Proof of Item 1:} By Lemma \ref{lemma1}, we know that the coordinates of $\mathcal{G}(\mathcal{F}^{in}\backslash W^s_{loc}(\CC_1))$ can be written as:
\begin{equation}
\nonumber
\left\{
\begin{array}{l}
\bigbreak
\phi_1^{out}=\Xi(r_1^{in} \cos \phi_2^{in},r_1^{in} \sin \phi_2^{in}) -\xi\omega_1 \ln(1- r_1^{in}) \pmod{2\pi}\\
\bigbreak
r_2^{out}=1-(1-r_1^{in})^\delta\\
\bigbreak
\phi_2^{out}=\phi_2^{in}- \xi \omega_2\ln(1- r_1^{in}) \pmod{2\pi},
\end{array}
\right.
\end{equation}
where $1-\hat{\delta}\leq r_1^{in} < 1$ and $0\leq \phi_2^{in}\leq 2\pi$. Omitting  $\mod{2\pi}$, we can write $\phi_1^{out}$ as a function of $r_2^{out}$ and $ \phi_2^{out}$  as follows:
\begin{equation}
\label{extra}
\begin{array}{l}
\phi_1^{out} = \Upsilon(r_2^{out}, \phi_2^{out}) = \Xi\left[ u (r_2^{out}, \phi_2^{out}), v (r_2^{out}, \phi_2^{out}) \right] - \frac{\xi\omega_1}{\delta}\ln(1- r_2^{out}).\\
\end{array}
\end{equation}
\bigbreak
\bigbreak
\noindent\textbf{Claim 1:} The following property holds:
\begin{equation}
\nonumber
\frac{\partial \Upsilon}{\partial r_2^{out}}(r_2^{out}, \phi_2^{out}) =  \frac{1}{1-r_2^{out}} \left[\frac{\xi \omega_1}{\delta} +R(r_2^{out},\phi_2^{out})\right],
\end{equation}
where $\dpt \lim_{r_2 \to 1^{-}}R(r_2^{out},\phi_2^{out})=0$. The convergence to $0$ is uniform with respect to $\phi_2^{out}$.\\

\begin{proof}
From \eqref{extra} it follows that
\begin{equation*}
\frac{\partial \Upsilon}{\partial r_2^{out}}(r_2^{out}, \phi_2^{out})
= \frac{\partial \Xi}{\partial u}(u, v)\frac{\partial u}{\partial r_2^{out}}(r_2^{out}, \phi_2^{out})
+ \frac{\partial \Xi}{\partial v}(u, v)\frac{\partial v}{\partial r_2^{out}}(r_2^{out}, \phi_2^{out})
+ \frac{1}{1-r_2^{out}}  \frac{\xi \omega_2}{\delta}.
\end{equation*}
Taking into account the expressions for $u$ and $v$,  one gets:

\begin{eqnarray}
\frac{\partial u}{\partial r_2^{out}}(r_2^{out}, \phi_2^{out})
&=&
\frac{dM}{d r_2^{out}} \cos \mathcal{W}
- M \sin \mathcal W \, \frac{\partial \mathcal W}{\partial r_2^{out}} \nonumber
\\
&=&
\frac{1}{1-r_2^{out}}
\left[ \frac{1}{\delta} (1-r_2^{out})^{1/\delta} \cos \mathcal{W}
+ \frac{\xi\omega_2}{\delta} (1-(1-r_2^{out})^\frac{1}{\delta}) \sin \mathcal{W} \right]\nonumber
\end{eqnarray}
and
\begin{eqnarray}
\frac{\partial v}{\partial r_2^{out}}(r_2^{out}, \phi_2^{out})
&=&
\frac{dM}{d r_2^{out}} \sin \mathcal{W}
+ M \cos \mathcal W \, \frac{\partial \mathcal W}{\partial r_2^{out}} \nonumber
\\
&=&
\frac{1}{1-r_2^{out}}
\left[ \frac{1}{\delta} (1-r_2^{out})^{1/\delta} \sin \mathcal{W}
- \frac{\xi\omega_2}{\delta} (1-(1-r_2^{out})^\frac{1}{\delta}) \cos \mathcal{W} \right]. \nonumber
\end{eqnarray}
Therefore
\begin{equation}
\nonumber
\frac{\partial \Upsilon}{\partial r_2^{out}}(r_2^{out}, \phi_2^{out}) =  \frac{1}{1-r_2^{out}} \left[\frac{\xi \omega_1}{\delta} +R(r_2^{out},\phi_2^{out})\right],
\end{equation}
where $R$ can be written as
\begin{equation*}
R(r_2^{out},\phi_2^{out})=R_1(r_2^{out},\phi_2^{out})+R_2(r_2^{out},\phi_2^{out})
\end{equation*}
with
\begin{equation*}
R_1(r_2^{out},\phi_2^{out})
=\frac{1}{\delta} (1-r_2^{out})^{1/\delta}
\left[ \frac{\partial \Xi}{\partial u}(u, v) \cos \mathcal{W}
+ \frac{\partial \Xi}{\partial v}(u, v) \sin \mathcal{W}  \right]
\end{equation*}
and
\begin{equation*}
R_2(r_2^{out},\phi_2^{out})
=\frac{\xi \omega_2}{\delta} (1-(1-r_2^{out})^{1/\delta})
\left[ \frac{\partial \Xi}{\partial u}(u, v) \sin \mathcal{W}
- \frac{\partial \Xi}{\partial v}(u, v) \cos \mathcal{W}  \right].
\end{equation*}
Since $\Xi$ is $C^1$, it easily follows that $\dpt \lim_{r_2^{out} \to 1^{-}} R_1(r_2^{out},\phi_2^{out})=0$. On the other hand,  
\begin{equation}
\frac{\partial \Xi}{\partial u}(u, v) \sin \mathcal{W}
- \frac{\partial \Xi}{\partial v}(u, v) \cos \mathcal{W}
=
\frac{1}{M(r_2^{out})}\left[ \frac{\partial \Xi}{\partial u}(u, v) v
- \frac{\partial \Xi}{\partial v}(u, v) u  \right].
\label{eq:aux1}
\end{equation}
Since the $C^1$--map $\Xi$ is constant when $u^2+v^2=1$  (by hypothesis), it follows that
\begin{equation*}
\left[ \frac{\partial \Xi}{\partial u}(u, v) v
- \frac{\partial \Xi}{\partial v}(u, v) u  \right]\equiv 0
\end{equation*}
for $(u,v)$ on the unit circle. Taking into account that $\dpt \lim_{r_2^{out} \to 1^-} M(r_2^{out})=1$ we conclude that the expression in \eqref{eq:aux1} tends to $0$ as $r_2^{out}\to 1^-$ and hence
$\dpt \lim_{r_2^{out} \to 1^{-}} R_2(r_2^{out},\phi_2^{out})=0$. This proves Claim~1. \\
\end{proof}

\noindent \textbf{Claim 2:} The following property holds:
\begin{equation}
\nonumber
\lim_{r_2^{out} \to 1^-}\frac{\partial \Upsilon}{\partial \phi_2^{out}}(r_2^{out}, \phi_2^{out}) = 0.
\end{equation}
The convergence to $0$ is uniform  with respect to $\phi_2^{out}$.

\begin{proof}
We use similar arguments to those of the proof of Claim 1. First, from \eqref{extra} it follows that
\begin{equation*}
\frac{\partial \Upsilon}{\partial \phi_2^{out}}(r_2^{out}, \phi_2^{out})
= \frac{\partial \Xi}{\partial u}(u, v)\frac{\partial u}{\partial \phi_2^{out}}(r_2^{out}, \phi_2^{out})
+ \frac{\partial \Xi}{\partial v}(u, v)\frac{\partial v}{\partial \phi_2^{out}}(r_2^{out}, \phi_2^{out})
\end{equation*}
Now, taking into account the expressions for $u$ and $v$ we easily get
\begin{equation*}
\frac{\partial \Upsilon}{\partial \phi_2^{out}}(r_2^{out}, \phi_2^{out})
= -\frac{\partial \Xi}{\partial u}(u, v)v(r_2^{out}, \phi_2^{out})
+ \frac{\partial \Xi}{\partial v}(u, v)u(r_2^{out}, \phi_2^{out}).
\end{equation*}
Since $\Xi$ is a $C^1$--map which is constant on the unit circle and $\|(u,v)\|\to 1$ as $r_2^{out} \to 1^-$, then we get the result.
\end{proof}

\bigbreak

In particular, if $1-r_2^{out}$ is small enough, namely if there exists $\hat\delta_1\gtrsim 0$ such that  $1-\hat\delta<1-\hat\delta_1 <r_2^{out}\leq 1$, then $\frac{\partial \Upsilon}{\partial r_2^{out}}(r_2^{out}, \phi_2^{out})>0$ and hence $ \Upsilon$ is injective as function of $r_2^{out}$. \\

For $\delta_2>0$ small, define the two cylindrical annulus in $\Sigma_2^{out}$:
$$
C_1=\{(\phi_1^{out}, r_2^{out}, \phi_2^{out}) \in \Sigma_2^{out}:
\quad \phi_1^{out}, \phi_2^{out} \in \RR\pmod{2\pi}
\quad \text{and} \quad 1-\delta_2<r_2^{out}\leq  1\}
$$
and
 \begin{eqnarray*}
C_2&=&\biggl\{(\phi_1^{out}, r_2^{out}, \phi_2^{out}) \in \Sigma_2^{out}:
\quad \phi_1^{out}, \phi_2^{out} \in \RR\pmod{2\pi}
\\
&& \text{and} \quad 1-\exp\left(\frac{\delta}{\xi \omega_2} \Upsilon (r_2^{out}, \phi_2^{out}) \right)
\leq r_2^{out}\leq 1 \biggl\}.
\end{eqnarray*}
 
\bigbreak

We would like to show that there exists a diffeomorphism between $C_1$ and $C_2$ transforming $\mathcal{G}(\mathcal{F}^{in}\backslash W^s_{loc}(\CC_1))$ into a reference spiralling sheet as in  Harterich \cite{Hart}. Define:
$$
\eta(\phi_1^{out}, r_2^{out}, \phi_2^{out})=(\tilde{\phi}_1^{out}, \tilde{r}_2^{out}, \tilde{\phi}_2^{out})
$$
with $(\phi_1^{out}, r_2^{out}, \phi_2^{out})\in C_1$ and
\begin{equation}
\nonumber
\left\{
\begin{array}{l}
\bigbreak
\tilde{\phi}_1^{out}= {\phi}_2^{out}\\
\bigbreak
\tilde{r}_2^{out}= h(r_2^{out}, \phi_2^{out}) \\
\bigbreak
\tilde{\phi}_2^{out}={\phi}_1^{out}\\
\end{array}
\right.
\end{equation}
where\\
\begin{equation}
\nonumber
h(r_2^{out}, \phi_2^{out})=
\left\{
\begin{array}{ll}
\dpt 1-\exp \left(-\frac{\delta}{\xi \omega_1}\Upsilon(r_2^{out}, \phi_2^{out})\right) &\qquad  \text{if} \qquad r_2^{out} \neq 1\\
1&\qquad  \text{if} \qquad r_2^{out} = 1.
\end{array}
\right.
\end{equation}
\bigbreak
\bigbreak
\noindent\textbf{Claim 3:} The map $\eta$ is injective.
\begin{proof}
Assume that $\eta(\phi_1^{out}, r_2^{out}, \phi_2^{out})=\eta(\hat \phi_1^{out}, \hat r_2^{out}, \hat \phi_2^{out})$. Therefore: \\
\begin{equation}
\nonumber
\left\{
\begin{array}{l}
\bigbreak
\phi_2^{out}=\hat \phi_2^{out}\\
\bigbreak
h(r_2^{out}, \phi_2^{out})=h(\hat r_2^{out}, \hat \phi_2^{out})=h(\hat r_2^{out}, \phi_2^{out})\\
\bigbreak
\phi_1^{out}=\hat \phi_1^{out}.\\
\end{array}
\right.
\end{equation}
Since $\Upsilon$ is injective for fixed $\phi_2^{out}$ and sufficiently small $r_2^{out}$, the same property holds for $h$, as we can assume that $\delta_2$ is also sufficiently small. Consequently, we have $r_2^{out}=\hat r_2^{out}$, leading to the conclusion that $\eta$ is injective. \\
\end{proof}

\noindent\textbf{Claim 4:} The sets $C_1$ to $C_2$ are diffeomorphic by  the map $\eta$. 

\begin{proof}
One only needs to prove that $h$ is a $C^1$--function. From Claim 1, it easily follows that
\begin{eqnarray}
\nonumber \frac{\partial h}{\partial r_2^{out}}(r_2^{out}, \phi_2^{out})
&=&\frac{\delta}{\xi \omega_1}
\frac{\partial \Upsilon}{\partial r_2^{out}}(r_2^{out}, \phi_2^{out})
\exp \left(-\frac{\delta}{\xi \omega_1} \Upsilon(r_2^{out}, \phi_2^{out})\right) \\
\nonumber &=&
\frac{1}{1-r_2^{out}} \left[1 + \frac{\delta}{\xi \omega_1} R(r_2^{out},\phi_2^{out})\right]
\exp \left(-\frac{\delta}{\xi \omega_1}\Upsilon(r_2^{out}, \phi_2^{out})\right).
\end{eqnarray}
Now, from the definition of $\Upsilon$ in \eqref{extra}, we get
\begin{eqnarray}
\nonumber
\frac{\partial h}{\partial r_2^{out}}(r_2^{out}, \phi_2^{out})
&=&
\frac{1}{1-r_2^{out}} \left[1 + \frac{\delta}{\xi \omega_1} R(r_2^{out},\phi_2^{out})\right]
\exp \left(
-\frac{\delta}{\xi \omega_1}\Xi(u,v)
+\ln(1-r_2^{out})
\right)
\\
\nonumber
&=&
\left[1 + \frac{\delta}{\xi \omega_1} R(r_2^{out},\phi_2^{out})\right]
\exp \left( -\frac{\delta}{\xi \omega_1}\Xi(u,v) \right).
\end{eqnarray}
On the other hand,
\begin{eqnarray}
\nonumber \frac{\partial h}{\partial \phi_2^{out}}(r_2^{out}, \phi_2^{out})
&=&\frac{\delta}{\xi \omega_1}
\frac{\partial \Upsilon}{\partial \phi_2^{out}}(r_2^{out}, \phi_2^{out})
\exp \left(-\frac{\delta}{\xi \omega_1} \Upsilon(r_2^{out}, \phi_2^{out})\right).
\end{eqnarray}
It follows from Claim 1 and Claim 2, respectively, that
$$
\lim_{r_2^{out} \rightarrow 1^-}\frac{\partial h}{\partial r_2^{out}}(r_2^{out}, \phi_2^{out})=	\exp \left(-\frac{\delta C}{\xi \omega_1}  \right)  \neq 0,
$$
where $C=\Xi(u,v)$ when $\|(u,v)\|=1$, and 
$$
\lim_{r_2^{out}\rightarrow 1}\frac{\partial h}{\partial \phi_2^{out}}(r_2^{out},\phi_2^{out})=0.
$$
The convergence to $0$ is uniform with respect to $\phi_2^{out}$. Therefore $\eta$ is a $C^1$ map. 
On the other hand,
\begin{equation*}
D\eta(\phi_1^{out}, r_2^{out}, \phi_2^{out})=
\left(\begin{array}{ccc}
0 & 0 &  1  \\
\\
0 & \frac{\partial h}{\partial r_2^{out}}(r_2^{out},\phi_2^{out})  
  & \frac{\partial h}{\partial \phi_2^{out}}(r_2^{out},\phi_2^{out})  \\
\\
1&0&0
\end{array}\right)
\end{equation*}
and hence,
\begin{equation*}
\lim_{r_2^{out} \to 1^-} D\eta(\phi_1^{out}, r_2^{out}, \phi_2^{out})=
\left(\begin{array}{ccc}
0 & 0 &  1  \\
\\
0 & \exp \left( -\frac{\delta C}{\xi \omega_1}  \right)  &   0  \\
\\
1 & 0 & 0
\end{array}\right).
\end{equation*}
Therefore, $||D\eta||$ is bounded away from zero in $C_1$ and, consequently, $\eta$ is a $C^1$--diffeomorphism. 
 \end{proof}

The set $C_1$ plays the role of $V$ in Definition \ref{ss_def2}. Finally, note that, by (\ref{extra}), we have:
$$
\tilde{r}_2^{out}= 1-\exp \left(-\frac{\delta}{\xi \omega_1}\Upsilon(r_2^{out}, \phi_2^{out})\right) = 1-\exp \left(-\frac{\delta}{\xi \omega_1} \phi_1^{out}\right)
= 1-\exp \left(-\frac{\delta}{\xi \omega_1} \tilde{\phi}_2^{out}\right).
$$
Since, in polar coordinates, the curve defined by $(\tilde{r}_2^{out}, {\tilde\phi}_2^{out})$ defines a spiral $\mathcal{S}$ accumulating on a curve lying on $W^u_{loc}(\CC_2) \cap \Sigma^{out}_2$, the proof of the first item ends.\\

\textit{Proof of Item 2:} The proof runs along the same lines as those of the previous item. \\

\textit{Proof of Item 3:} For $i\in \{1,2\}$, we know that:
$$S_i= \mathcal{G}^{-1}(C_i^{out}\backslash W^u(\CC_2))$$
where \\
$$
 C_i^{out}=\left\{ ( \phi_1^{out}, r_2^{out}, \phi_2^{out}):
r_2^{out}\in [1-\varepsilon^{out},1],
\phi_1^{out} \in [0,2\pi], \phi_2^{out} \in [\theta_i^{out}-\varepsilon^{out},\theta_i^{out}+\varepsilon^{out}] \right\} \subset \Sigma_2^{out}.
$$
\bigbreak
Consider the following two disks limiting  $C_i^{out}$:
$$
 \mathcal{D}_1 =\left\{ ( \phi_1^{out}, r_2^{out}, \phi_2^{out}):
r_2^{out}\in [1-\varepsilon^{out},1],
\phi_1^{out} \in [0,2\pi], \phi_2^{out} = \theta_i^{out}-\varepsilon^{out}] \right\} \subset \Sigma_2^{out}.
$$
and 
$$
 \mathcal{D}_2 =\left\{ ( \phi_1^{out}, r_2^{out}, \phi_2^{out}):
r_2^{out}\in [1-\varepsilon^{out},1],
\phi_1^{out} \in [0,2\pi], \phi_2^{out} = \theta_i^{out}+\varepsilon^{out}] \right\} \subset \Sigma_2^{out}.
$$

\bigbreak
 Each of these disks limiting $C_i^{out}$, is sent by $\mathcal{G}^{-1}$ into a spiralling sheet accumulating on $W^{s}_{loc}(\CC_1) \cap \Sigma_1^{in}$. Hence $S_i$ is a set limited by spiralling sheets inside the hollow cylinder $\Sigma_1^{in}$. Then, by definition, it is a scroll contained in $\Sigma_1^{in}$ accumulating on $W^{s}_{loc}(\CC_1) \cap \Sigma_1^{in}$.
 \end{proof}

\begin{figure}
\begin{center}
\includegraphics[height=6.9cm]{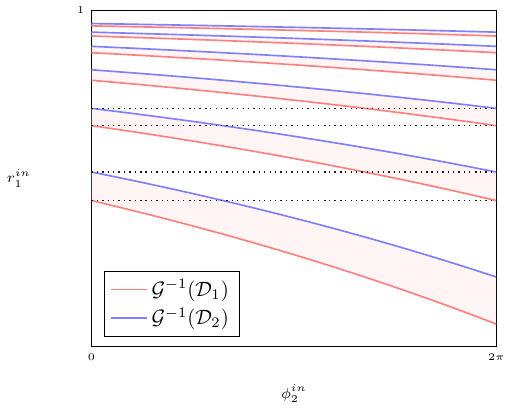}
\end{center}
\caption{\small { Sketch of the intersection between the scroll $\mathcal{G}^{-1}(C_1^{out}\setminus W^{u}(\mathbf{C}_2))$ and a section with $\phi_1^{in}$ constant. The sets $\mathcal{D}_1$ and $\mathcal{D}_2$ are those used in the proof of Item 3 in Proposition \ref{My_Prop5}.}}
\label{scroll}
\end{figure}
 Figure \ref{scroll} illustrates Item 3 in Proposition \ref{My_Prop5}. The intersection between the scroll $\mathcal{G}^{-1}(C_1^{out}\setminus W^{u}(\mathbf{C}_2))$ and a section in $\Sigma_1^{in}$ (red shaded region) is shown. It is limited by $\mathcal{G}^{-1}(\mathcal{D}_1)$ (red lines) and $\mathcal{G}^{-1}(\mathcal{D}_2)$ (blue lines), where $\mathcal{D}_i$, with $i=1,2$, are the sets defined in the proof of Item 3.
Now we are in conditions to prove Proposition \ref{infinite number of connections}.

\subsection{Proof of Proposition \ref{infinite number of connections}}
\label{ss:proof1}
 For $\gamma > 0$ (small), by Hypothesis \textbf{(P6)}, the set $W^u(\CC_2)\cap \Sigma_1^{in}$ defines a surface intersecting transversally the set $W^s(\CC_1)\cap \Sigma_1^{in}$. Then, making use of Proposition \ref{My_Prop5}, the set $\mathcal{G}(W^u(\CC_2)\cap \Sigma_1^{in})$ contains two spiralling sheets accumulating on $W^u(\CC_2)\cap \Sigma_2^{out}$. In particular, it intersects the set $W^s(\CC_1)\cap \Sigma_2^{out}$  along infinitely many curves. Each curve corresponds to a two-dimensional heteroclinic connection from $\CC_2$ to $\CC_1$.

\section{The tridimensional horseshoe}
\label{s:horseshoe}
In this section, we recall the main steps of the construction of the invariant horseshoe given by Wiggins \cite{Wiggins} adapted to our purposes. We address the reader to Section 2.3 of the latter book  for more details in the definitions.

\subsection{Preliminaries for the construction}
\label{slices and slabs}
\begin{figure}
\begin{center}
\includegraphics[height=6.9cm]{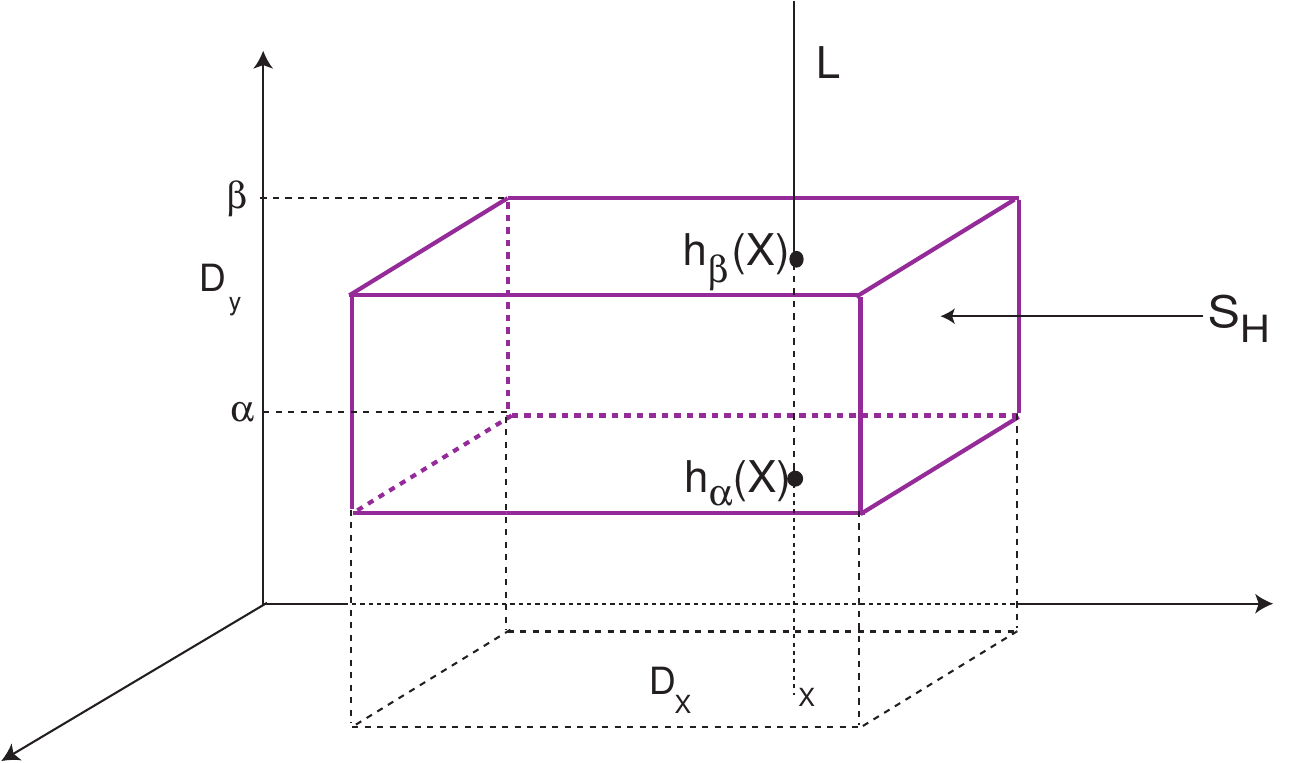}
\end{center}
\caption{\small Illustration of the sets $D_X$, $D_y$, $ {S_H}$ (horizontal slab) and Property \textbf{(P)}.}
\label{box}
\end{figure}

Let $D\subset \RR^3$ be a compact and connected 3-dimensional set of $\RR^3$. 
Define $D_X$ and $D_y$ the projection of $D$ onto $\RR^2$ and $\RR$ as:
$$D_X=\{X\in \RR^2: \text{ for which there exists } y\in \RR \text{ with } (X, y)\in D\}\subset \RR^2$$
and
$$D_y=\{y\in \RR: \text{ for which there exists } X\in \RR^2 \text{ with } (X, y)\in D\}\subset \RR. $$
In the case illustrated in Figure \ref{box},  $D_X$ is a closed and connected two-dimensional rectangle contained in $\RR^2$ and $D_y$ is a compact interval of $\RR$. In what follows, $I$ denotes a set of indices.

\begin{definition}
A $\mu_h$-horizontal slice $\mathcal{H}$ is a set defined as    $$\mathcal{H}= \{(X, h(X))\in \RR^3: X \in D_X\} \subset D,$$ where  $h: D_X \rightarrow \RR$ is a map for which 
   there exists $\mu_h\in \RR^+_0$ such that  $$|h(X_1)-h(X_2)|\leq \mu_h \|X_1-X_2\|, \quad \text{for all} \quad X_1, X_2\in D_X.$$
 
\bigbreak
A $\mu_v$-vertical slice $\mathcal{V}$ is a set defined as 
$$ \mathcal{V}= \{(v(y), y)\in \RR^3: y \in I_y\} \subset D,$$
where $v: D_y \rightarrow \RR^2$ is a map for which  there exists $\mu_v\in \RR^+_0$ such that: $$\|v(y_1)-v(y_2)\|\leq \mu_v |y_1-y_2|,  \quad \text{for all} \quad y_1, y_2\in D_y.$$ 
 \end{definition}
\bigbreak
 
For $\mu_h>0$ fixed, let $\mathcal{H}\subset D$ be a $\mu_h$-horizontal slice and let $J\subset D$ homeomorphic to a closed interval, intersecting $\mathcal{H}$ at any (but only one) point. Let $\mathcal{H}^\alpha$, $\alpha\in I$, be the set of horizontal slices that intersect the boundary of $J$, say $\partial J$, and have the same domain of $\mathcal{H}$. Now, consider the following set:
$$
S_\mathcal{H} = \{ (X,y) \in \RR^2\times \RR: X \in D_X \text{ and } y \text{ has the Property \textbf{(P)} }\}
$$
where, for a set of one-dimensional indices $I$:\\
\medbreak
 \textbf{Property (P): } For each $X\in D_X$, given any line $L$ through $(X,y)$ with L parallel to vertical line $X=0$, then $L$ intersects the points $(X, h_\alpha (X))$ and $(X, h_\beta (X))$ for some $\alpha, \beta \in I$ with $(X, y)$ between these two points along $L$.  The sets $h_\alpha (X)$ and $h_\beta (X)$ should be seen as horizontal slices.\\

\begin{definition}
A $\mu_h$\emph{-horizontal slab} is defined to be the topological closure of $S_\mathcal{H}$.  \\
 \end{definition}
 
From now on, without loss of generality, we assume that $S_\mathcal{H}$ is closed. The \emph{vertical boundary} of a $\mu_h$-horizontal slab $S_\mathcal{H}$, denoted by $\partial_v S_\mathcal{H}$, is defined as 
$$\partial_v S_\mathcal{H} = \{(X,y)\in S_H: X\in \partial D_X\}.$$

The \emph{horizontal boundary} of a $\mu_h$-horizontal slab $S_\mathcal{H}$, denoted by $\partial_h S_\mathcal{H}$, is defined by $\partial_h S_\mathcal{H} = \partial S_\mathcal{H} \backslash \partial_v S_\mathcal{H}.$ \\

\begin{definition}
Let $S_\mathcal{H} ^1$ and $S_\mathcal{H} ^2$ be two $\mu_h$-horizontal slabs. We say that $S_\mathcal{H} ^1$ intersects $S_\mathcal{H} ^2$ \emph{fully} if $S_\mathcal{H} ^1\subset S_\mathcal{H} ^2$ and  $\partial_v S_\mathcal{H} ^1\subset \partial_v S_\mathcal{H} ^2$.\\
\end{definition}
\bigbreak
\begin{definition}
Fix $\mu_v > 0$. Let $S_\mathcal{H}$ be a $\mu_h$-horizontal slab and let  $\mathcal{V}$ be a $\mu_v$ vertical slice contained in $S_\mathcal{H}$ such that
$$
\partial  { \mathcal{V}}\subset \partial_h S_\mathcal{H}.
$$
Let $J\subset \mathcal{H}$ be a 2-dimensional disk intersecting {  $\mathcal{V}$} at any (but just one) point of { $\mathcal{V}$}, and let {  $\mathcal{V}^\alpha$}, $\alpha \in I$, be the set of all $\mu_v$-vertical slices that intersect the boundary of $J$ with $\partial \mathcal {V}^\alpha \subset \partial _h S_\mathcal{H}$.   A $\mu_v$-vertical slab is defined to be the topological closure of the set 
\begin{eqnarray*}
\lefteqn{S_\mathcal{V}= \lbrace (X,y)\in \RR^3:}
\\
&& (X,y) \text{  is contained in the interior of the set bounded by {  ${V}^\alpha$}, with $\alpha\in I$, and $\partial_h S_\mathcal{H}$} \rbrace.
\end{eqnarray*}
\end{definition}
\noindent As before, we assume that $S_\mathcal{V}$ is closed.\\

\begin{definition}
Let $S_\mathcal{V}$ be a $\mu_v$--vertical slab. The horizontal boundary of $S_\mathcal{V}$, denoted $\partial_h S_\mathcal{V}$, is defined to be $S_\mathcal{V} \cap \partial_h S_\mathcal{H}$. The vertical boundary of $S_\mathcal{V}$, denoted $\partial_v S_\mathcal{V}$, is defined to be $\partial S_\mathcal{V} \backslash \partial_h S_\mathcal{V}$.\\
\end{definition}

\begin{definition}
\label{def_width}
The width of a $\mu_h$-horizontal slab $S_\mathcal{H}$, denoted by $d(S_\mathcal{H})$, is defined as:
\begin{equation}
\label{dist}
d(S_\mathcal{H})= \sup_{X\in D_X, \alpha, \beta \in I} |h_\alpha (X)-h_\beta(X)|.
\end{equation}
In a similar way, the  width of a $\mu_v$-vertical slab $S_\mathcal{V}$ is 
\begin{equation}
\label{dist}
d(S_\mathcal{V})= \sup_{y\in D_y, \alpha, \beta \in I} |v_\alpha (y)-v_\beta(y)|.
\end{equation}
 \end{definition}

\subsection{Proof of Theorem \ref{Th2}}

\label{suspended horseshoes}
\begin{figure}
\begin{center}
\includegraphics[height=8.0cm]{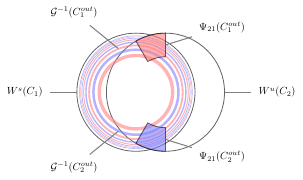}
\end{center}
\caption{\small Intersection of the scrolls $\mathcal{G}^{-1}(C^{out}_1\setminus W_{loc}^{u}(\CC_2))$ (red annulus) and $\mathcal{G}^{-1}(C^{out}_2\setminus W_{loc}^{u}(\CC_2))$ (blue annulus) with a section where $\phi_1^{in}$ remains constant within $\Sigma_1^{in}$. The images of $C_i^{out}$ (for $i=1,2$) under $\Psi_{21}$ on $\Sigma_1^{in}$ are also depicted.  The red and blue annuli correspond to the images under $\mathcal{R}^{-1}_{\gamma}$ of $\Psi_{21}(C_1^{out})$ and $\Psi_{21}(C_2^{out})$, respectively. Inside the intersection sets $\mathcal{G}^{-1}(C^{out}_i\setminus W_{loc}^{u}(\CC_2)) \cap \Psi_{21}(C_i^{out})$, $i=1,2$, horseshoes are formed. }
\label{scroll_side_view}
\end{figure}

For $i\in \{1,2\}$, for $\gamma \gtrsim 0$ fixed, in this section we focus our attention on the dynamics of 
\begin{equation}
\label{def_S}
S_i= \mathcal{G}^{-1}(C^{out}_i\setminus W_{loc}^{u}(\CC_2))  \subset \Sigma_{1}^{in}
\end{equation}
 defined in Subsection \ref{ss:notation1}.  Geometrically, as a consequence of Proposition \ref{My_Prop5}, the set $S_i$ is a scroll (see Figure \ref{scroll_side_view}). 
  To simplify the readers' task we revisit few results,  whose arguments will be required in the sequel. We prove the existence of  invariant sets in $$S=S_1\cup S_2,$$ accumulating on $W_{loc}^s(\CC_1) \cap W^u(\CC_2) \cap \Sigma_1^{in}$, where $\mathcal{R}_\gamma|_S$ is topologically conjugate to a shift under (at least) two symbols -- we remind the readers that $\mathcal{R}_\gamma$ stands for the first return map to $\Sigma_1^{in}$ (cf. Subsection \ref{ss:notation2}). The construction is based on the generalized Conley-Moser conditions \cite{Wiggins}, which provide sufficient conditions for the existence of invariant sets where the dynamics is conjugated to a full shift. 

For $i\in \{1,2\}$, let $C_i^{out}\subset \Sigma_2^{out}$ be a portion of the hollow cylinder  with radius $1-\varepsilon^{out} \leq r_2^{out} \leq 1$  given by (\ref{ciout}). Given  $\varepsilon^{out}>0$ fixed, for $N \in \NN$  large enough  (see Remark \ref{N large 1} below) and fixed, 
we define the set:
\begin{eqnarray*}
M_{N,1}^{out}=\left\{ (\phi_1^{out}, r_2^{out}, \phi_2^{out}):
r_2^{out}\in [1-a_{N},1-a_{N+1}],
\phi_1^{out} \in [0,2\pi], \phi_2^{out} \in [-\varepsilon^{out},\varepsilon^{out}] \right\} \subset C^{out}_1,
\end{eqnarray*}
and
\begin{eqnarray*}
M_{N,2}^{out}=\left\{ (\phi_1^{out}, r_2^{out}, \phi_2^{out}):
r_2^{out}\in [1-a_{N},1-a_{N+1}],
\phi_1^{out} \in [0,2\pi], \phi_2^{out} \in [\pi -\varepsilon^{out},\pi +\varepsilon^{out}] \right\} \subset C^{out}_2,
\end{eqnarray*}

where
\begin{equation}
\label{a_n}
a_{N}=\exp\left(  - \frac{2\pi N\delta }{\xi\omega_2}\right).
\end{equation}
It is easy to check that $\dpt \lim_{N \to +\infty}(1-  a_N)=1$ and for all $N\in \NN$ (large enough), $a_N>a_{N+1}$.

\begin{remark}
\label{N large 1}
We precise what is meant by ``$N\in \NN$ large enough''. 
In order that the sets $M_{N,1}^{out}$ and $M_{N,2}^{out}$ are well defined in $C_1^{out}$ and $C_2^{out}$, respectively,  then $N\in \NN$ must be such that:
 
\begin{equation*}
1-a_{N}>1-\varepsilon^{out} \quad
\Leftrightarrow \quad \exp\left(  - \frac{2\pi N \delta }{\xi\omega_2}\right)<\varepsilon^{out} \quad \Leftrightarrow \quad N>-\ln (\varepsilon^{out}) \frac{\xi\omega_2}{2\pi \delta}>0.
\end{equation*}
 
\end{remark}

\begin{figure}
\begin{center}
\includegraphics[height=5.0cm]{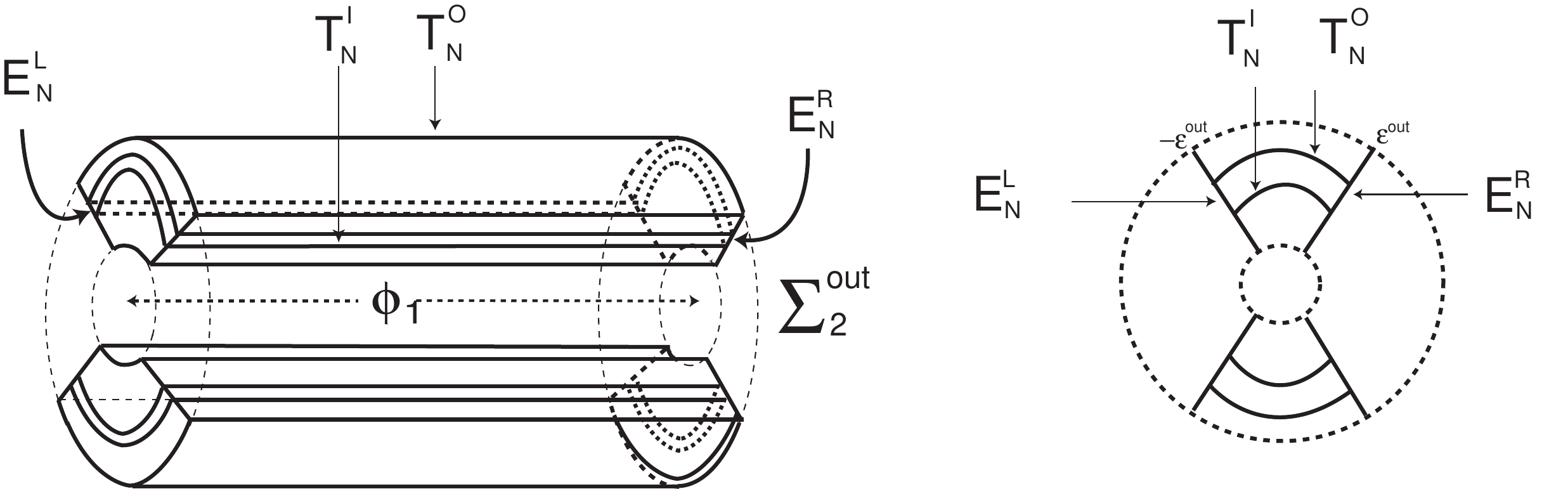}
\end{center}
\caption{\small  Topological boundary of $M_{N,1}^{out}\subset \Sigma_2^{out}$. Assuming $M_{N,1}^{out} \equiv M_{N}^{out}\subset \Sigma_2^{out}$, we have: $\partial M_{N}^{out} = E_N^{L} \cup E_N^{R}  \cup T_N^{I} \cup T_N^{O}$, where the letters $L$, $R$, $I$, $O$ mean Left, Right, Inner and Outer, respectively.}
\label{annulus}
\end{figure}

In order to simplify the notation, we perform the following analysis for  $M_{N}^{out}\equiv M_{N,1}^{out}$. Computations for   $M_{N,2}^{out}$ are completely similar with suitable changes in the $\phi_2^{out}$ coordinate (instead of $\phi_2^{out} =\pm \varepsilon^{out}$, we should consider $\phi_2^{out} =\pi\pm \varepsilon^{out}$).
 As shown in Figure \ref{annulus}, the topological boundary of $M_{N}^{out}\subset \Sigma_2^{out}$, denoted by $\partial M_{N}^{out}$, can be written as
$$
\partial M_{N}^{out} = E_N^{L} \cup E_N^{R}  \cup T_N^{I} \cup T_N^{O},
$$
(the letters $L$, $R$, $I$, $O$ mean Left, Right, Inner and Outer, respectively) with: \\
\begin{eqnarray*}
E_N^{L}&=&\left\{ (\phi_1^{out},r_2^{out} \phi_2^{out})\in \Sigma_2^{out}:
r_2^{out}\in [1-a_{N},1-a_{N+1}],
\phi_1^{out} \in [0,2\pi], \phi_2^{out} =-\varepsilon^{out} \right\},
\\
E_N^{R}&=&\left\{ (\phi_1^{out},r_2^{out} \phi_2^{out})\in \Sigma_2^{out}:
r_2^{out}\in  [1-a_{N},1-a_{N+1}],
\phi_1^{out} \in [0,2\pi], \phi_2^{out}= \varepsilon^{out} \right\},
\\
T_N^{I}&=&\left\{ (\phi_1^{out},r_2^{out} \phi_2^{out})\in \Sigma_2^{out}:
r_2^{out}=1-a_{N},
\phi_1^{out} \in [0,2\pi], \phi_2^{out} \in [-\varepsilon^{out},\varepsilon^{out}] \right\} \quad \text{and }
\\
T_N^{O}&=&\left\{ (\phi_1^{out},r_2^{out} \phi_2^{out})\in \Sigma_2^{out}:
r_2^{out}=1-a_{N+1},
\phi_1^{out} \in [0,2\pi], \phi_2^{out} \in [-\varepsilon^{out},\varepsilon^{out}] \right\}.
\end{eqnarray*}
\bigbreak
 According to the definitions given in Subsection \ref{slices and slabs}, it is easy to check the set $M_{N}^{out}$ is part of a horizontal slab across $\Sigma_2^{out}$. The vertical (resp. horizontal) boundaries of $M_{N}^{out}$ are defined by $E_N^L$ and $E_N^R$ (resp. $T_N^I$ and $T_N^O$). The surfaces $T_N^I$ and $T_N^O$ may be defined as graphs of functions $r_2^{out}=h(\phi_1^{out},\phi_2^{out})$ where $h$ is  the constant map equal to $1-a_N$ and $1-a_{N+1}$, respectively. \\

 For $i\in \{1,2\}$, define:
\begin{equation}
\label{sn}
 \mathcal{S}_{N,i}= \mathcal{G}^{-1}(M_{N,i}^{out})\subset \Sigma_1^{in}.
 \end{equation}
 
  Using the coordinates of the different components of $\partial M_{N}^{out}$ we prove the following technical and useful result:

\begin{lemma} 
\label{lemma_imageofsolidannulus}
The following equalities hold:
\begin{eqnarray*}
{ \mathrm{(i)}}\:\mathcal{G}^{-1}(E_N^{L})= & \{ (r_1^{in},\phi_1^{in}, \phi_2^{in}):
&       r_1^{in}= 1-(1-r_2^{out})^\frac{1}{\delta},
\\
& &     \phi_1^{in} = \phi_1^{out}+\frac{\xi \omega_1}{\delta}\ln (1-r_2^{out}),
 \\
& &     \phi_2^{in}  =-\varepsilon^{out} + \frac{\xi \omega_2}{\delta}\ln (1-r_2^{out}),
 \\
& &   \text{where}\quad   r_2^{out}\in [1-a_N, 1-a_{N+1}], \phi_1^{out} \in [0,2\pi] \} \subset \Sigma_1^{in},
\\
\\
{  \mathrm{(ii)}}\:\mathcal{G}^{-1}(E_N^{R})= & \{ (r_1^{in},\phi_1^{in}, \phi_2^{in}):
&       r_1^{in}= 1-(1-r_2^{out})^\frac{1}{\delta},
\\
& &     \phi_1^{in} = \phi_1^{out}+\frac{\xi \omega_1}{\delta}\ln (1-r_2^{out}),
 \\
& &     \phi_2^{in}  =\varepsilon^{out} + \frac{\xi \omega_2}{\delta}\ln (1-r_2^{out}),
 \\
& &    \text{where}\quad  r_2^{out}\in [1-a_N, 1-a_{N+1}], \phi_1^{out} \in [0,2\pi] \} \subset \Sigma_1^{in},\
\\
\\
{  \mathrm{(iii)}}\:\mathcal{G}^{-1}(T_N^{I})= &  \{ ( r_1^{in},  \phi_1^{in},\phi_2^{in}):
&       r_1^{in}=1- \exp\left( -\frac{ 2\pi N}{\xi \omega_2}\right),
\\
& &     \phi_1^{in} = \phi_1^{out}-\frac{\omega_1}{\omega_2}2\pi N,
\\
& &     \phi_2^{in} = \phi_2^{out}-2\pi N,
\\
& &     {  \text{where}\quad   \phi_1^{out} \in [0, 2\pi], \phi_2^{out}\in [-\varepsilon^{out},\varepsilon^{out}]\} \subset \Sigma_1^{in},}
\\
\\
{  \mathrm{(iv)}}\:\mathcal{G}^{-1}(T_O^{I})= &  \{ ( r_1^{in},  \phi_1^{in},\phi_2^{in}):
&       r_1^{in}=1- \exp\left( -\frac{ 2\pi (N+1)}{\xi \omega_2}\right),
\\
& &     \phi_1^{in} = { \phi_1^{out}}-\frac{\omega_1}{\omega_2}2\pi (N+1),
\\
& &     \phi_2^{in} = \phi_2^{out}-2\pi (N+1),
\\
& &     {  \text{where}\quad   \phi_1^{out} \in [0, 2\pi], \phi_2^{out}\in [-\varepsilon^{out},\varepsilon^{out}]\}
 \subset \Sigma_1^{in},}
 \\
\end{eqnarray*}
\end{lemma}

\begin{figure}
\begin{center}
\includegraphics[height=8.9cm]{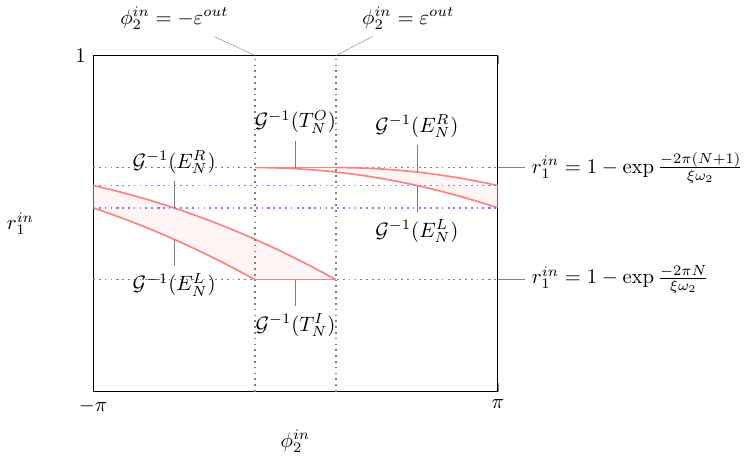}
\end{center}
\caption{\small Illustration of Lemma \ref{lemma_imageofsolidannulus}. Projection in the plane $(\phi_2^{in}, r_1^{in})$ of $\Sigma_1^{in}$ of the set $ \mathcal{S}_{N,i}= \mathcal{G}^{-1}(M_{N,i}^{out})\subset \Sigma_1^{in}$, for $i=1,2$. We have emphasised the role of the pre-images of its boundaries drawn in Figure \ref{annulus}. Vertical sides are identified.}
\label{side_view_new}
\end{figure}

\begin{proof}
We concentrate our attention in the first and third items. The proof of the others runs along the same lines with appropriate changes.  We suggest the reader to follow the proof by observing Figures \ref{annulus} and \ref{side_view_new}.\\

(i) 
Using the expression \eqref{return2} of Lemma \ref{lemma1} and the fact that $(\phi_1^{out},r_2^{out}, \phi_2^{out})  \in E_N^L$, one knows that:
\begin{eqnarray*}
 r_1^{in}&=& 1-(1-r_2^{out})^{\frac{1}{\delta}} \\
 \phi_1^{in} &=& \phi_1^{out}   +\frac{\xi \omega_1}{\delta}\ln (1-r_2^{out}) \pmod{2\pi} \\
 \phi_2^{in} &=& -\varepsilon^{out}   +\frac{\xi \omega_2}{\delta}\ln (1-r_2^{out}) \pmod{2\pi}
  \end{eqnarray*}

(ii) Similar to (i) by substituting $\phi_2^{out}$ by $\varepsilon^{out}$.\\

(iii) Using the expression \eqref{return2} and the fact that   $(\phi_1^{out},r_2^{out}, \phi_2^{out})\in T_N^I$, one knows that:

\begin{eqnarray*}
 r_1^{in}&=& 1-(1-r_2^{out})^{\frac{1}{\delta}} = 1-a_N^{\frac{1}{\delta}}
 = 1-\exp\left(-\frac{2\pi N}{\xi \omega_2}\right)\\
 \phi_1^{in} &=& \phi_1^{out}   +\frac{\xi \omega_1}{\delta}\ln a_N= \phi_1^{out}   -\frac{ \omega_1}{\omega_2}2\pi N  \pmod{2\pi} \\
 \phi_2^{in} &=& \phi_2^{out}   +\frac{\xi \omega_2}{\delta}\ln a_N =  \phi_2^{out}   -2\pi N \pmod{2\pi}
\end{eqnarray*}

(iv) Similar to (iii) but replacing $N$ by $N+1$.
\end{proof}

In the following result, we prove that the image under $\mathcal{G}^{-1}$ of $M_N^{out}$ covers $\Sigma_1^{in}$ (at least once) with respect to both angular coordinates. We assume that the set $[0,2\pi]$ is endowed with the usual quotient topology of $\RR / 2\pi$.
\begin{lemma}
\label{bijection1}
 The following assertions hold: 
\begin{enumerate}
\item The components $\phi_1^{in}$ and $\phi_2^{in}$ of $ \mathcal{G}^{-1}(E_N^{L})$ are surjective in $[0,2\pi]$.
\item The components $\phi_1^{in}$ and $\phi_2^{in}$ of $ \mathcal{G}^{-1}(E_N^{R})$ are surjective in $[0,2\pi]$.
\item The components $\phi_1^{in}$ and $\phi_2^{in}$ of $ \mathcal{G}^{-1}(T_N^{I})$ are surjective in $[0,2\pi]$ and $[-\varepsilon^{out}, \varepsilon^{out}]$, respectively.
\item The components $\phi_1^{in}$ and $\phi_2^{in}$ of $ \mathcal{G}^{-1}(T_N^{O})$ are surjective in $[0,2\pi]$ and $[-\varepsilon^{out}, \varepsilon^{out}]$, respectively.
\end{enumerate}
\end{lemma}

\begin{proof}
Again, we concentrate our attention in the first and third items. The proof of the others runs along the same lines with appropriate changes. \\
\begin{enumerate}
\item Using Lemma \ref{lemma_imageofsolidannulus} and the fact that $(\phi_1^{out},r_2^{out}, \phi_2^{out})\in E_N^L$, one knows that:
\begin{eqnarray*}
  \phi_1^{in} &=& \phi_1^{out}   +\frac{\xi \omega_1}{\delta}\ln (1-r_2^{out})   \\ 
\end{eqnarray*}
For a fixed value of $r_2^{out} \in [1-a_N, 1-a_{N+1}]$, since $\phi_1^{out}\in [0,2\pi]$, we get that $  \phi_1^{in}$ is surjective in $[0,2\pi]$.
 On the other hand, as function of $r_2^{out}$, $\phi_2^{in}$ is a decreasing function. Since $\phi_2^{in}=-\varepsilon^{out}-2\pi N$ for $r_2^{out}=1-a_N$ and $\phi_2^{in}=-\varepsilon^{out}-2\pi (N+1)$ for $r_2^{out}=1-a_{N+1}$, it follows that $\phi_2^{in}$ is surjective over $[0,2\pi]$.\\
%

\item Similar to (1) by replacing $\phi_2^{out}$ by $\varepsilon^{out}$. \\

\item
%
  Using Lemma \ref{lemma_imageofsolidannulus}
\begin{eqnarray*}
\phi_1^{in} &=& \phi_1^{out}-\frac{\omega_1}{\omega_2} 2 \pi N  \\
\phi_2^{in} &=& \phi_2^{out}-2\pi N
\end{eqnarray*}
with $\phi_1^{out}\in [-\varepsilon^{out},\varepsilon^{out}]$ and $\phi_2^{out}\in [0,2\pi]$,
forcing that $\phi_1^{in}$ is surjective in $[-\varepsilon^{out}, \varepsilon^{out}]$ and $\phi_2^{in}$ is surjective in $[0,2\pi]$. \\
 
 \item Similar to (iii).
\end{enumerate}
\end{proof}
 
\bigbreak
By Proposition \ref{My_Prop5}, for each $i\in \{1,2\}$,  the set $S_i= \mathcal{G}^{-1}(  C_i^{out}\setminus W^u(\CC_2))\subset \Sigma_1^{in}$ is a scroll accumulating on $W^s_{loc}(\CC_1)\cap \Sigma_1^{in}$. 
For   $i\in \{1,2\}$, the family of sets $\mathcal{S}_{N,i}\subset \Sigma_1^{in}$ (see \eqref{sn}) provides an infinite collection of pieces in the interior of the scroll $S_i$ accumulating (with $N$) on $W^s_{loc}(\CC_1)\cap \Sigma_1^{in}$.

For $N\in \NN$ fixed (large enough), note that $\mathcal{S}_{N,i}$ is limited by the lift of a two-dimensional  annulus contained in $\Sigma_1^{in}$. More precisely
$$
\mathcal{S}_{N,i} \subset \left\{ (r_1^{in},\phi_1^{in}, \phi_2^{in}) \in \Sigma_1^{in}:
1-b_{N}\leq r_1^{in}\leq 1-b_{N+1} \right\} \quad \text{where} \quad b_{N}=\exp\left( -\frac{  2\pi N}{\xi \omega_2}\right).
$$
 
 Observe that there exists a diference between $\mathcal{S}_{N,i}$ and $S_i$, where $N\in \NN$  and  $i\in \{1,2\}$ -- the second is a scroll accumulating on $W^s(\CC_1)\cap \Sigma_1^{in}$ and the first is a proper subset of the second. 
 \\

\begin{figure}
\begin{center}
\includegraphics[height=6.2cm]{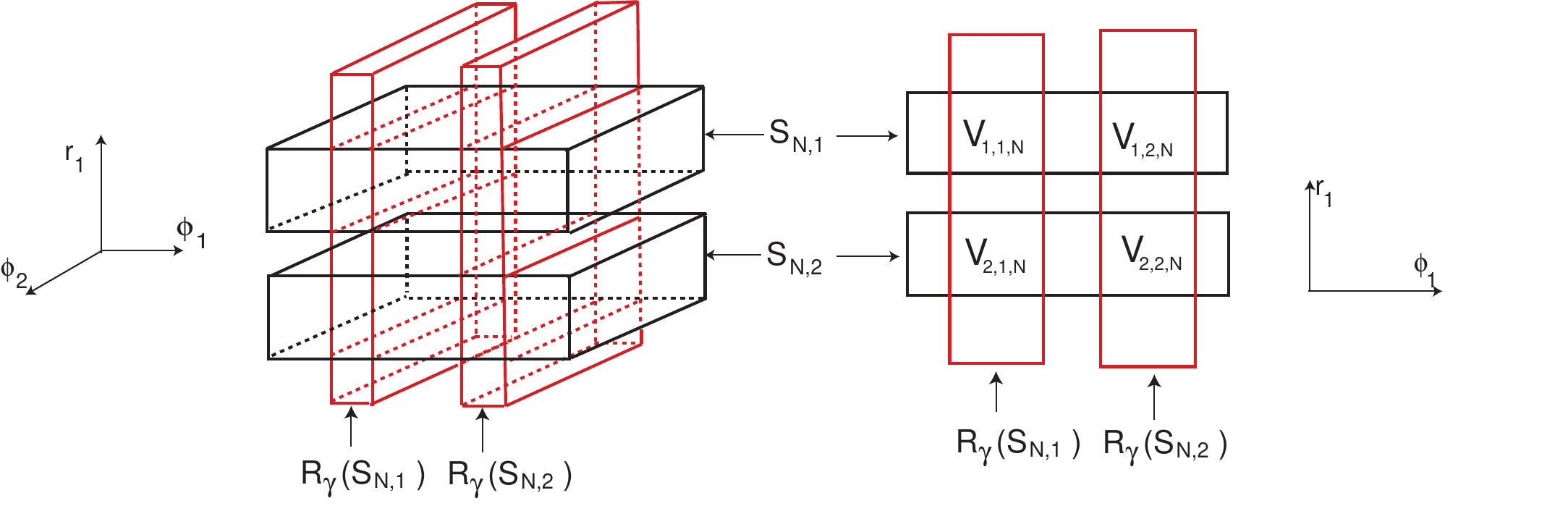}
\end{center}
\caption{\small Representation of the tridimensional horseshoe and its projection into the coordinates $(r_1^{in}, \phi_2^{in})$ of $\Sigma_1^{in}$. The red sets represents $\mathcal{R}_\gamma(\mathcal{S}_{N,2})$. Note that $V_{j,i, N}:=\mathcal{R}_\gamma(\mathcal{S}_{N,i})\cap\mathcal{S}_{N,j}$ may be seen as a vertical slab. Note that the sets $\mathcal{S}_{N,1}$ and $\mathcal{S}_{N,2}$ have been flattened.} 
\label{horseshoe1}
\end{figure}

 For $i, j  \in \{1,2\}$  and $N\in \NN$ fixed and large enough, let us consider the set depicted in Figure \ref{horseshoe1}: $$V_{j,i, N}:=\mathcal{R}_\gamma(\mathcal{S}_{N,i})\cap\mathcal{S}_{N,j}\subset \Sigma_1^{in},$$ where $\mathcal{R}_\gamma=\Psi_{2,1}\circ \mathcal{G}$ is the return map introduced in   Subsection \ref{global transition}.  From the previous construction, as a consequence of  Hypothesis \textbf{(P6)},  for all  $N$ large, we get $V_{j,i,N} \neq \varnothing$. Moreover, each $V_{j,i,  N}$ has  one connected component.    These sets are fully intersecting vertical slabs across $\mathcal{S}_{N,j}$. \\
 Define now $$H_{i, j, N}=\mathcal{R}_\gamma^{-1}(V_{j,i, N})$$ 
 
 \begin{remark} The following information may be useful to the reader.
 \begin{enumerate}
 \item  Since $V_{j,i, N}:=\mathcal{R}_\gamma(\mathcal{S}_{N,i})\cap\mathcal{S}_{N,j}$, then it may be useful the following equality:  $$H_{i, j, N} =\mathcal{S}_{N,i}\cap \mathcal{R}_\gamma^{-1}(\mathcal{S}_{N,j}).$$
 \item For each $N\in \NN$, the first subscript in $H_{i,j, N}$ informs which particular horizontal slab the set is in; the second indicates the vertical slab  the set is mapped under $\mathcal{R}_\gamma^{-1}$.
 \end{enumerate}
 \end{remark}

 The next result claims that $H_{i, j, N}$ is a fully intersecting horizontal slab across $\mathcal{S}_{N,i}$.

\begin{lemma} [Slab Condition]
\label{lemma12} For each  $i,j\in \mathbb{N}$ large enough, the set $H_{i,j, N}$
is a fully intersecting horizontal slab in $\mathcal{S}_{N,i}$.
\end{lemma}

\begin{proof}
Let us fix (once for all) the set $V_{j, i, N}$. One knows that
$$
\Psi_{21}^{-1}(V_{j,i, N})=\left\{ (\phi_1^{out},r_2^{out} \phi_2^{out}):
r_2^{out}\in [1-a_{N},1-a_{N+1}],
\phi_1^{out} \in [0,2\pi],   H_L \leq \phi_2^{out} \leq H_R \right\},
$$
where $H_R, H_L:D \to \RR$ are $C^1$ maps defined in $$D =\{(x,y)\in \RR^2: 1-a_{N}\leq x^2+y^2 \leq 1-a_{N+1}\}$$  such that $H_L<H_R$.
The sets 
$$
\hat{E}^L =\left\{ (\phi_1^{out},r_2^{out} \phi_2^{out}):
r_2^{out}\in [1-a_{N},1-a_{N+1}],
\phi_1^{out} \in [0,2\pi],    \phi_2^{out} = H_L \right\},
$$
and
 $$
\hat{E}^R =\left\{ (\phi_1^{out},r_2^{out} \phi_2^{out}):
r_2^{out}\in [1-a_{N},1-a_{N+1}],
\phi_1^{out} \in [0,2\pi],    \phi_2^{out} = H_R \right\}
$$
define the vertical boundaries of $\Psi_{21}^{-1}(V_{j,i,N})$. Combining Lemmas \ref{lemma_imageofsolidannulus} and \ref{bijection1} (first and second items), this implies that the sets $\mathcal{G}^{-1}(\hat{E}^L)$ and $\mathcal{G}^{-1}(\hat{E}^R)$ are surfaces  cutting transversally in $\mathcal{S}_{N,i}$. Hence, by construction, the set $H_{i,j, N}$
is a fully intersecting horizontal slab in $\mathcal{S}_{N,i}$.
\end{proof}
 
The image under $\mathcal{R}_\gamma$ of the two horizontal boundaries of each $H_{i,j,N}$ is a horizontal boundary of $V_{j,i, N}$. Moreover,  the image of each vertical boundary of $H_{i,j,N}$ is a  vertical boundary of $\mathcal{R}_\gamma(H_{i,j,N})$. Then:
 
 \begin{lemma}
 $\partial_v V_{j, i, N} \subset \partial_v \mathcal{R}_\gamma(H_{i, j, N})$ and $\mathcal{R}_\gamma^{-1}(\partial_v   V_{j, i, N} ) \subset \partial_v H_{i, j, N}$.
 \end{lemma}

We have checked the topological part of the Conley-Moser conditions, namely, all properties included in Hypothesis (A1) of \cite[Theorem 2.3.3]{Wiggins}. Now we need to obtain estimations of the rates of contraction and expansion of $H_{i,j,N}$ under $\mathcal{R}_\gamma$ along the horizontal and vertical directions. These analytical conditions require uniform (but not constant) contraction in the horizontal directions and expansion in the vertical direction. In what follows, we denote by $d(H)$ (resp. $d(V)$) the width of a horizontal (resp. vertical) slab $H$ (resp. $V$) -- see Definition \ref{def_width} above.  

\begin{lemma} [Hyperbolicity condition] 
\label{hyp_condition}
The following assertions hold for $i, j\in \{1,2\}$: \\
\begin{enumerate}
\item If $H$ is a $\mu_h$-horizontal slab intersecting $\mathcal{S}_{N,i}$ fully, then $\tilde{\mathcal{S}}_{N,j}:= \mathcal{R}_\gamma^{-1}(H)\cap \mathcal{S}_{N,j}$ is a $\nu_h$-horizontal slab intersecting $\mathcal{S}_{N,j}$ fully and $d( \mathcal{\tilde{S}}_{N,j})\leq \nu_h d(H)$, for some $\nu_h \in\, \, ]\, 0, 1\,[$.\\ 
\medbreak
\item If $V$ is a $\mu_v$-vertical slab contained in $\mathcal{S}_{N,j}$ such that $V\subset V_{j,i,  N}$, then $\tilde{V}_{N,i}=\mathcal{R}_\gamma(V)\cap \mathcal{S}_{N, i}$ is a  $\nu_v$-vertical slab contained in $\mathcal{S}_{N, i}$ and $d( \tilde{V}_{N,i})<\nu_v d(V)$, for some $\nu_v \in\, \, ]\, 0, 1\,[$. \\
\end{enumerate}
\end{lemma}

\begin{proof}

We are going to prove item (1); the proof of the other is analogous by considering the expression \eqref{return2}. 
Given a horizontal slab $H$ fully intersecting $\mathcal{S}_{N,i}$ for some $i=1,2$, then
$$\mathcal{R}_\gamma^{-1}(H)=\bigcup_{ i=1,2}\bigcup_{j=1,2} \mathcal{R}_\gamma^{-1} (H \cap V_{j,i,N}).
$$
We will prove that there exists $\nu_h \in\,  ]0,1[$ such that $d(\mathcal{R}_\gamma^{-1} (H \cap V_{j,i,N}))<\nu_h  d(H)$. Note that the set $\Psi_{2,1}^{-1}(H \cap V_{j,i,N})$ is a subset in $C_i^{out}$ (by construction) which we can write as:
$$
\begin{array}{rcl}
\Psi_{21}^{-1}(H \cap V_{j,i,N})&=&\left\{ (r_2^{out}, \phi_2^{out},  \phi_1^{out}):r_2^{out}\in [1-a_{N},1-a_{N+1}], \phi_1^{out} \in [0,2\pi],\right.
\\
\\
&&
\left. H_L  \leq \phi_2^{out} \leq  H_R \right\} \subset M_{i,N}^{out}
\end{array}
$$
for some for some $C^1$--maps $H_L$ and $H_R$  defined on the annulus $$D=\{(u,v)\in\RR^2:\,\:\,1-a_N \leq u^2+v^2 \leq 1-a_{N+1}\}.$$ The maps $H_L$ and $H_R$ are bounded and, in order to simplify the next computations and without loss of generality, we assume that they are constant maps, namely, $H_L \equiv \tau_L$ and $H_R \equiv \tau_R$ with $\tau_L < \tau_R$. 
As in the proof of Lemma \ref{lemma12}, the sets 
$$
\hat{E}^L =\left\{ (\phi_1^{out},r_2^{out} \phi_2^{out}):
r_2^{out}\in [1-a_{N},1-a_{N+1}],
\phi_1^{out} \in [0,2\pi],    \phi_2^{out} = \tau_L \right\},
$$
and
 $$
\hat{E}^R =\left\{ (\phi_1^{out},r_2^{out} \phi_2^{out}):
r_2^{out}\in [1-a_{N},1-a_{N+1}],
\phi_1^{out} \in [0,2\pi],    \phi_2^{out} = \tau_R \right\}
$$
define the vertical boundaries of $\Psi_{21}^{-1}(H\cap V_{i,j,N})$. 
Omitting $\mod 2\pi$, using the expression \eqref{return2}, the horizontal boundary of $\mathcal{G}^{-1}\circ \Psi_{21}^{-1}(H\cap V_{i,j,N})(\hat{E}^R )$ satisfies the condition:
\begin{eqnarray}
\label{expressao1.1}
\phi_2^{in}&=& \tau_R  +\frac{\xi \omega_2}{\delta}\ln (1-r_2^{out}) \Leftrightarrow \\
\nonumber r_2^{out} &=& 1- \exp \left( \frac{(\phi_2^{in}-\tau_R)\delta}{\xi \omega_2}\right)
\end{eqnarray}
and then, using again \eqref{return2}, we get:
$$
r_1^{in}= 1- (1- r_2^{out})^\frac{1}{\delta}= 1- \exp \left( \frac{(\phi_2^{in}-\tau_R) }{\xi \omega_2}\right).
$$
A similar equality holds for  $\mathcal{G}^{-1}\circ \Psi_{21}^{-1}(H\cap V_{i,j,N})(\hat{E}^L )$. The height of  $\mathcal{R}_\gamma^{-1}(H\cap V_{i,j,N})$ is given by the difference of the $r_1^{in}$--coordinates. Indeed, 

\begin{eqnarray*}
  \exp \left( \frac{(\phi_2^{in}-\tau_R) }{\xi \omega_2}\right)-  \exp \left( \frac{(\phi_2^{in}-\tau_L) }{\xi \omega_2}\right) &=&     {\exp \left( \frac{\phi_2^{in} }{\xi \omega_2}\right)}\left[  \exp \left( \frac{ -\tau_R }{\xi \omega_2} \right)-  \exp \left( \frac{ -\tau_L }{\xi \omega_2} \right) \right] \\
  &\overset{\text{Lagrange}}=&      {\exp \left( \frac{\phi_2^{in} }{\xi \omega_2}\right)}\left[  \exp \left( \frac{ -c }{\xi \omega_2} \right)  \right] (\tau_R-\tau_L)\\
  &=&  \exp \left( \frac{(\phi_2^{in}-c) }{\xi \omega_2}\right)(\tau_R-\tau_L) \\
    &\overset{\eqref{expressao1.1}}=&   \exp \left( \frac{-c }{\xi \omega_2}\right)  \exp \left(  \tau_R+\frac{\xi \omega_2}{\delta} \ln (1-r_2^{out}) \right)(\tau_R-\tau_L) \\
    &\overset{\eqref{expressao1.1}}=& \exp \left( \frac{\tau_R-c}{\xi \omega_2}\right)   \left(\frac{\xi \omega_2}{\delta} \ln (1-r_2^{out}) \right)(\tau_R-\tau_L) \\
    &\leq &   {K}(1-r_2^{out})^\frac{\xi \omega_2}{\delta}  (\tau_R-\tau_L) \\
    &\overset{\eqref{a_n}}<& {K} a_{N}^\frac{\xi \omega_2}{\delta}  (\tau_R-\tau_L) \\ 
    \end{eqnarray*}
 where $c\in [\tau_L,\tau_R] \subset [\theta_i^{out}-\varepsilon^{out},\theta_i^{out}+\varepsilon^{out}]$   and $K\in \RR^+$. Since $a_N$ may be taken as small as required taking $N\in \NN$ large enough (cf. \eqref{a_n}), the result follows.

 \end{proof}
Lemma \ref{hyp_condition} shows that the condition (A2) of \cite[pp. 118]{Wiggins} is valid. Thus, by  \cite[Theorem 2.3.3]{Wiggins}, for each $N\in \NN$ (large enough) there exists an invariant set
$$
{\Lambda_\gamma}=\bigcap_{N \in \mathbb{N}} \bigcup_{i=-N}^{N}\mathcal{R}_\gamma^{-i}\left( \mathcal{S}_{N,1} \cup  \mathcal{S}_{N,2} \right) \subset \Sigma_1^{in}
$$
on which the first return map $\mathcal{R}_\gamma|_{\Lambda_\gamma}$ is topologically conjugate to a full shift over an alphabet of two symbols $\{1,2\}$.
 By construction, the maximal invariant set in $\mathcal{V}$ is the topological closure of $ \Lambda_\gamma  $ and it is a Cantor set with zero Lebesgue measure  because of Lemma \ref{hyp_condition}. It corresponds to the maximal invariant set by construction.
 Theorem \ref{Th2} is now proved. \\
 
 \begin{remark}
 The zero Lebesgue measure of ${\Lambda_\gamma}$ is consistent with Bowen Theorem which states that $C^2$--hyperbolic sets are either Anosov or have zero Lebesgue measure.
 \end{remark}

\begin{remark}
There are two moments where we  need ``$N$ large enough''. The first has been discussed in Remark \ref{N large 1}; the other has been used at the end of the proof of Lemma \ref{hyp_condition}. So,  Theorem \ref{Th2} holds for  $n$ larger than the maximum of these two thresholds. 
\end{remark}
\begin{remark}
The  proof of Theorem \ref{Th2}  allows us to conclude that when  $\gamma> 0$ (small) and $N,M\in \NN$ are large enough, the hyperbolic horseshoes $\Lambda_{N}$,  $\Lambda_{M}$ are heteroclinically related. More precisely, the unstable manifolds of the periodic orbits in $\Lambda_{N}$, are long enough to intersect the stable manifolds of the periodic points of $\Lambda_{M}$. That is, given two slabs, there exist periodic solutions jumping from one slab to another, and so the heteroclinic classes associated to the infinitely many horseshoes are not disjoint.   Regarding this subject, the chapter about the \emph{whiskers} of the horseshoes in   \cite{GST} is worthwhile reading.
  \end{remark}
 \section{Concluding remark}
  \label{s:concluding}
  
  Motivated by a certain type of unfolding of a Hopf-Hopf singularity, the so called \emph{Gaspard-type unfolding} defined in Subsection \ref{ss:Gaspard-type}, in this paper we have considered a one-parameter family $(f_\gamma)_{\gamma\geq0}$ of $C^3$--vector fields in $\RR^4$ whose flows exhibit a heteroclinic cycle associated to two  periodic solutions and a bifocus, all of them hyperbolic. The periodic solutions have complex (non-real) Floquet multipliers.   
   For $\gamma>0$, we have proved the existence of infinitely many linked horseshoes with zero Lebesgue measure. This is the content of the main result of this paper (Theorem \ref{Th2}).     By the way the construction has been performed  in Section \ref{s:horseshoe},  the hyperbolic Cantor set $\Lambda_\gamma$ has one contracting and two expanding directions. Heteroclinic switching appears as a consequence of the main result.
   
   The proof has been based in the rigorous checking of Conley and Moser conditions established in \cite{Wiggins}, using the techniques of \cite{ALR, IbRo, KR2010}. The novelty   lies on the existence of a periodic solution in dimension 4, which makes the description much more involving. Due to the transverse intersection of $W^u(\CC_2)$ and $W^s(\CC_1)$, for $\gamma>0$, the horseshoes whose existence is guaranteed by Theorem \ref{Th2} are persistent and does not depend on the saddle value of the nodes as the three-dimensional homoclinic orbit to a saddle-focus \cite{Shilnikov65}.    In Section \ref{ss:Gaspard-type}, we have shown the existence of tridimensional horseshoes in particular analytic unfoldings of a Hopf-Hopf singularity. In a different context, this technique has   been used in Section 10 of \cite{Rodrigues2022} to prove the existence of dynamics associated to Torus-breakdown  near a Hopf-zero singularity of Gaspard-type.

   The next natural problem to be solved is the existence of persistent strange attractors  (observable chaos in the sense defined by \cite{WO2011}) near the network $\Gamma_\gamma$. At this moment, assuming $\delta>1$, we have two conjectures: \\
   
   \textbf{Conjecture A:} Under the conditions of Theorem \ref{Th2}, there exist a sequence $(\gamma_i)_{i\in \NN}$ for which the vector field $\dot{x}=f_{\gamma_i}(x)$ exhibits a   heteroclinic tangency between $W^u(\CC_2)$ and $W^s(\CC_1)$. This would be 
   the scenario to guarantee the existence of strange attractors. In this case, we do not break the network $\Gamma_\gamma$ associated to  $\bold{O}$, $\CC_1$ and   $\CC_2$.
    Another way to prove the existence of strange attractors  and strange repellers (via tangencies in the first return map) is to follow the ideas of Homburg \cite{hom}:  moving genericaly the saddle-value of the network is enough to obtain tangencies. We defer this task for future work. \\

   \textbf{Conjecture B:} Under the conditions of Theorem \ref{Th2}, if $f_\gamma$ satisfies Hypotheses \textbf{(P1)}--\textbf{(P7)}  then there exists a diffeomorphism $G$ arbitrarily $C^1$-close to $\mathcal{R}_\gamma$, exhibiting (H\'enon-type) strange attractors  and/or infinitely many sinks.  In this case, we guess that it is possible to find a generalised Tatjer tangency which prompts the existence of a map close to $\mathcal{R}_\gamma$, where strange attractors may be found. The itinerary of this proof would follow the arguments of \cite{Rodrigues2020}.\\

Although we are not able to prove for the moment the existence of two-dimensional strange attractors close to $\Gamma_\gamma$, the scenario described by a vector field satisfying \textbf{(P1)}--\textbf{(P7)} is the natural setting in which topological two-dimensional strange attractors might occur. Therefore, it would be a challenge to find  analytically two-dimensional strange attractors when the network is broken. We defer these tasks for future work.


\begin{thebibliography}{99}

\bibitem{ACL}
Aguiar, M.A.D., Castro, S.B.D., Labouriau, I.S.,  \emph{Dynamics near a heteroclinic network},   Nonlinearity {\bf 18},
391--414, 2005.

\bibitem{ALR}
 Aguiar, M.A.D., Labouriau, I.S., Rodrigues, A.A.P.,  \emph{Switching near a heteroclinic network of rotating nodes}, 
Dyn. Syst {\bf 25}(1), 75--95, 2010.



\bibitem{BIS} Baldom\'a, I.,  Ib\'a\~nez, S., Seara, T.,  \emph{Hopf-Zero singularities truly unfold chaos},  Commun. Nonlinear Sci. Numer. Simul. \textbf{84} 105162, 2020.
 
\bibitem{Barrientos}
  Barrientos, P. G., Ib\'a\~nez, S.,  Rodr\'iguez, J.A., \emph{Heteroclinic cycles arising in generic unfoldings
of nilpotent singularities},  J. Dyn. Diff. Equations {\bf 23},  999--1028, 2011.

\bibitem{Barrientos_book}
Barrientos, P. G., Ib\'a\~nez, S., Rodrigues, A. A., Rodr\'iguez, J. \'A, \emph{ Emergence of Chaotic Dynamics from Singularities}., 32th Brazilian Mathematics Colloquium, Instituto Nacional de Matem\'atica Pura e Aplicada (IMPA), Rio de Janeiro, 2019.

\bibitem{Belitski}
  Belitskii, G.R., \emph{Functional equations, and conjugacy of local diffeomorphisms of finite smoothness class}, 
Funkcional. Anal. i Prilozen {\bf 7},  17--28, 1973.


\bibitem{BC91}  Benedicks, M.,   Carleson, L.,  \emph{The dynamics of the H\'enon map}. Ann. of Math. (2), \textbf{133}(1):73--169, 1991.


%

\bibitem{Birkhoff}
Birkhoff, G.D.,  \emph{ Nouvelles recherches sur les systmes dynamiques},  Memoriae Pont. Acad. Sci. Novi Lyncaei, \textbf{1} 85--216, 1935.
 
\bibitem{cdi2024}
Casas, P. S., Drubi, F.,  Ib\'a\~nez, S.,  \emph{Invariant manifolds in a reversible hamiltonian system: the tentacle-like geometry}. Commun. Nonlinear Sci. Numer. Simul., \textbf{138}, 108189, 2024.

\bibitem{CY08}
Campbell,  S.A., Yuan, Y.,  \emph{Zero singularities of codimension two and three in delay differential equations}. Nonlinearity, \textbf{21}(11): 267--2691, 2008.


\bibitem{din25}
Drubi, F., Ib\'a\~nez, S.,  Noriega, D., \emph{Hopf-Hopf bifurcations in a coupling of Fitzhugh-Nagumo systems}, Preprint, 2025. 

\bibitem{Drubi}
  Drubi, F., Ib\'a\~nez, S., Rivela, D., \emph{Chaotic behavior in the unfolding of Hopf-Bogdanov-Takens singularities},  Discrete Cont. Dyn. Syst. - B, \textbf{25}(2), 599-615, 2020.
 
 \bibitem{DIR07}
Drubi, F., Ib\'a\~nez, S.,  Rodr\'iguez, \emph{Coupling leads to chaos}. J. Diff. Equations, {\bf 239} (2): 371--385, 2007.
 
\bibitem{DIR11}
Drubi, F., Ib\'a\~nez, S.,  Rodr\'iguez, \emph{Hopf-pitchfork singularities in coupled systems}. Physica D: Nonlinear Phenomena, {\bf 240} (9-10), 825--840, 2011.

\bibitem{dumibakok2005}
  Dumortier, F., Ib\'a\~nez, S.,  Kokubu, H., \emph{Cocoon bifurcations in three-dimensional reversible vector fields},  {\it Nonlinearity} {\bf 19}  305--328, 2006.

\bibitem{dumibakok2001}
  Dumortier, F., Ib\'a\~nez, S., Kokubu, H., \emph{New aspects in the unfolding of the
nilpotent singularity of codimension three},  Dyn. Syst. {\bf 16}(1), 63--95, 2001.


\bibitem{dumibakoksim}
  Dumortier, F., Ib\'a\~nez, S., Kokubu, H., Sim\'o, C.,  \emph{About the unfolding of a Hopf-zero singularity}, 
Discrete Contin. Dyn. Syst. {\bf 33}, 4435--4471, 2013.


 

\bibitem{FS}
  Fowler, A.C.,  Sparrow, C.T., \emph{Bifocal homoclinic orbits in four dimensions},  {\it Nonlinearity,} {\bf 4}, 1159--1182, 1991.


\bibitem{Gaspard} Gaspard, P.,  \emph{ Local birth of homoclinic chaos},  Physica D: Nonlinear Phenomena, \textbf{62}(1-4)  94--122, 1993.

\bibitem{Glendinning}
  Glendinning, P., \emph{Differential equations with bifocal homoclinic orbits}, Internat. J. Bifur. Chaos.{\bf 7}(1), pp. 27--37, 1997.


\bibitem{GS1}
  Glendinning, P.,  Sparrow, C., \emph{Local and Global Behaviour near Homoclinic Orbits},  J. Stat. Phys. {\bf 35}, pp. 645--696, 1984.

 

\bibitem{GST} Gonchenko, S.V., Shilnikov, L.P., Turaev, D.V.,  \emph{Dynamical phenomena in systems with structurally unstable Poincar\'e homoclinic orbits}, Chaos 6, 15--31, 1996.



\bibitem{Hart}
  Harterich, J., \emph{Cascades of reversible homoclinic orbits to a bifocus equilibrium}, 
  Phys. D {\bf 112},  187--200, 1998.





\bibitem{hom}
  Homburg, A.J., \emph{Periodic attractors, strange attractors and hyperbolic dynamics near homoclinic orbits to saddle-focus equilibria}, Nonlinearity {\bf 15},  1029--1050, 2002.


\bibitem{HK}
 Homburg, A.J.,  Knobloch, J., \emph{Switching homoclinic networks}, 
 Dyn. Syst. {\bf 23},  351--358, 2010.


\bibitem{IbRo} Ib\'a\~nez, A., Rodrigues, A.A.P.;  \emph{On the dynamics near a homoclinic network to a bifocus:
switching and horseshoes}, Int. J. Bifurcation Chaos 2 {\bf 5}, 1530030, [19 pages], 2015.


\bibitem{IbanezRodriguez}
 Ib{\'a}{\~n}ez, S.,  Rodr{\'{\i}}guez, J.A., \emph{Shil'nikov configurations in any generic unfolding of the
nilpotent singularity of codimension three on $\mathbb{R}^3$,}  J. Diff. Equations {\bf 208},  147--175, 2005.


\bibitem{KR2010}
Kirk, V.,  Rucklidge, A., \emph{The effect of symmetry breaking on the dynamics near a structurally stable heteroclinic cycle between equilibria and a periodic orbit}, Dynamical Systems,  {\bf 23(1)} 43--74, 2008.

\bibitem{KL2009} Koltsova, O., Lerman, L.,  \emph{Hamiltonian dynamics near nontransverse homoclinic orbit to saddle-focus equilibrium}, Discrete Cont. Dyn. Sys., {\bf 25}, 3, 883--913, 2009.

  \bibitem{kuznetsov1998elements} Kuznetsov, Y. A., \emph{ Elements of applied bifurcation theory}, Springer, 1998.


\bibitem{LG}
 Laing, C., Glendinning, P., \emph{Bifocal homoclinic bifurcations},   Phys. D {\bf 102},  1--14, 1997.
 
 
\bibitem{LR2017}
 Labouriau, I.S., Rodrigues, A.A.P., \emph{On Takens' last problem: tangencies and time averages near heteroclinic networks}, Nonlinearity 30(5), 1876, 2017.

\bibitem{druetal2024}
Mayora-Cebollero, A., Jover-Galtier, J.A., Drubi, F., Ib\'a\~nez, S., Lozano, A., Mayora-Cebollero, C., Barrio, R., \emph{Almost synchronization phenomena in the two and three coupled Brusselator systems}, Physica D: Nonlinear Phenomena {\bf 472}, 134457, 2025.

\bibitem{MV93}
  Mora, L.,  Viana, M.,  \emph{Abundance of strange attractors},  Acta Math., \textbf{171}(1):1--71, 1993.

 

\bibitem{pm} Palis, J.,  de Melo, W., \emph{Geometric theory of dynamical systems: An introduction}, (Springer-Verlag, New York-Berlin), 1982.


\bibitem{Poincare}
 Poincar\'e, H., \emph{Sur le probleme ds trois corps et les \'eequations de la dynamique}, Acta Math., 13:1--270, 1890.
 
  
\bibitem{PR_book}
 Pumari\~no, A.,  Rodr\'iguez, J.A., \emph{Coexistence and persistence of strange attractors}, 
 Lecture Notes in Math., {\bf 1658}, (Springer, Berlin), 1997.


\bibitem{PR}
 Pumari\~no, A.,  Rodr\'iguez, J.A., \emph{Coexistence and persistence of infinitely many strange attractors}, 
Ergodic Theory Dynam. Systems {\bf 21}(5), 1511--1523, 2001.

\bibitem{Rodrigues3}
 Rodrigues, A., \emph{Repelling dynamics near a Bykov cycle},   J. Dyn. Diff. Equations {\bf 25}(3), 605--625, 2013.

\bibitem{Rodrigues2020}
Rodrigues, A., \emph{Strange attractors and wandering domains near a homoclinic cycle to a bifocus},  Journal of Diff. Equations, \textbf{269}(4), 3221--3258, 2020.

 \bibitem{Rodrigues2022}
 Rodrigues, A., \emph{Unfolding a Bykov attractor: from an attracting torus to strange attractors}. Journal of Dynamics and Differential Equations, 34(2), 1643-1677, 2022.

\bibitem{RLA}
  Rodrigues, A., Labouriau, I., Aguiar, M., \emph{Chaotic Double Cycling}, 
 Dyn. Syst. {\bf 26}(1), pp. 199--233, 2011.



\bibitem{Shilnikov65}
  Shilnikov, L.P., \emph{A case of the existence of a denumerable set of periodic motions}, 
 { Sov. Math. Dokl} {\bf 6}, 163--166, 1965.


\bibitem{Shilnikov67A}
  Shilnikov, L.P., \emph{The existence of a denumerable set of periodic motions in four dimensional space in an extended neighbourhood of a saddle-focus},   {Sov. Math. Dokl} {\bf 8}(1),  54--58, 1967.

\bibitem{Shilnikov70}
  Shilnikov, L.P., \emph{A contribution to the problem of the structure of an extended neighbourhood of a rough equilibrium state of saddle-focus type}, '  {\it Math. USSR Sb.} {\bf 10}(1), 91--102, 1970.

\bibitem{Shilnikov et al} L. Shilnikov, A. Shilnikov, D. Turaev, L. Chua, \emph{Methods of Qualitative Theory in Nonlinear Dynamics 1}, World Scientific Publishing Co., 1998.


\bibitem{Smale}
 Smale, S., \emph{Differentiable dynamical systems}. Bull. Amer. Math. Soc., 73:747-- 817, 1967.

\bibitem{Takens74}
 Takens, F.,  \emph{Singularities of vector fields}, Publications Math\'ematiques de l'IHES, 43(1):47--100, 1974.

\bibitem{Tresser}
  Tresser, C. , \emph{About some theorems by L. P. Shilnikov},  {Ann. Inst. H. Poincar\'e.} {\bf 40},  441--461, 1984.


\bibitem{WO2011} 
 Wang, Q., Ott, W., \emph{Dissipative homoclinic loops of two-dimensional maps and strange attractors with one
direction of instability}, Communications on Pure and Applied Mathematics 64(11), 1439--1496, 2011.

\bibitem{Wiggins}
  Wiggins, S., \emph{Global Bifurcations and Chaos. Analytical Methods}, {Applied Mathematical Sciences} {\bf 73}, Springer-Verlag, New York, 1988.

\end{thebibliography}
\end{document}